\documentclass{amsart}

\title{Rigid geometry on projective varieties}
\author{Benjamin McKay}
\address{University College Cork \\ National University of Ireland}
\date{\today} % \\ MSC Primary: 57R17, 51A10; Secondary: 14N05}
\thanks{Thanks to Indranil Biswas, Sorin Dumitrescu, 
Maciej Dunajski and J. M. Landsberg for helpful comments.}

\usepackage{ifthen}
\usepackage{euscript}
\usepackage{varioref}
\usepackage{longtable}
\usepackage[all]{xy}
\xyoption{web}
\usepackage{booktabs}

\newtheorem{theorem}{Theorem}
\newtheorem{corollary}{Corollary}
\newtheorem{lemma}{Lemma}
\newtheorem{proposition}{Proposition}
\theoremstyle{remark}
\newtheorem{conjecture}{Conjecture}
\newtheorem{definition}{Definition}
\newtheorem{example}{Example}
\newtheorem{remark}{Remark}

\newlength{\setBracketHeight}

\newcommand{\hook}{\ensuremath{\mathbin{ \hbox{\vrule height1.4pt%
        width4pt depth-1pt \vrule height4pt width0.4pt depth-1pt}}}}
\newcommand{\Lm}[2]{\ensuremath{\Lambda^{#1} \left ( {#2} \right )}}

\newcommand{\R}[1]{\ensuremath{\mathbb{R}^{#1}}}
\newcommand{\C}[1]{\ensuremath{\mathbb{C}^{#1}}}

\newcommand{\Proj}[2][]{\ensuremath{\mathbb{#1 P}^{#2}}}

\newcommand{\SL}[1]{\ensuremath{\operatorname{SL}\left({#1}\right)}}
\newcommand{\Sp}[1]{\ensuremath{\operatorname{Sp}\left({#1}\right)}}

\newcommand{\SO}[1]{\ensuremath{\operatorname{SO}\left({#1}\right)}}
\newcommand{\PSL}[1]{\ensuremath{\operatorname{PSL}\left({#1}\right)}}
\newcommand{\PO}[1]{\ensuremath{\operatorname{PO}\left({#1}\right)}}

\newcommand{\LieB}{\ensuremath{\mathfrak{b}}}
\newcommand{\LieG}{\ensuremath{\mathfrak{g}}}
\newcommand{\LieH}{\ensuremath{\mathfrak{h}}}

\newcommand{\LieP}{\ensuremath{\mathfrak{p}}}

\newcommand{\LieSL}[1]{\ensuremath{\mathfrak{sl}\left({#1}\right)}}
\newcommand{\Sym}[2]{\ensuremath{\operatorname{Sym}^{#1}\left({#2}\right)}}
\newcommand{\s}[1]{\LieSL{2,\C{}}_{#1}}

\newcommand{\nForms}[2]{\ensuremath{\Omega^{#1}\left({#2}\right)}}
\DeclareMathOperator{\Ad}{Ad}
\DeclareMathOperator{\ad}{ad}
\newcommand{\pout}{{\scriptscriptstyle{\times}}}
\newcommand{\pin}{{\scriptscriptstyle{\bullet}}}

\newcommand{\Grnull}[2]{\ensuremath{\operatorname{Gr}_{\operatorname{null}}\left({#1},{#2}\right)}}

\def\latticebody{\drop{}}

\begin{document}
\begin{abstract}
We prove rigidity of various types of holomorphic geometric
structures on smooth complex projective varieties.
\end{abstract}

\maketitle
\tableofcontents

\section{Introduction}
All manifolds and maps henceforth are assumed complex analytic, and
all Lie groups, algebras, etc. are complex. 
We will prove a collection of global rigidity theorems
for holomorphic geometric structures. As an example:

\begin{theorem}\label{theorem:RigidityForQuaternionicContact}
Suppose that $M$ is a smooth connected complex projective variety
of complex dimension 7, bearing a holomorphic quaternionic contact
structure. Then $M=C_3/P$, where $C_3=\Sp{6,\C{}}$  
and $P \subset C_3$ is a certain complex parabolic subgroup. 
The holomorphic quaternionic contact structure is the standard
flat holomorphic quaternionic contact structure on $C_3/P$ (defined in
section~\vref{section:QuaternionicContact}).
\end{theorem}

\section{The main theorem}
\subsection{Definitions required to state the main theorem}

\subsubsection{Definition of Cartan geometries}
\begin{definition}
If \(E \to M\) is a principal right \(G\)-bundle,
we will write the right \(G\)-action as \(r_g e = eg\), where $e\in E$
and $g\in G$.
\end{definition}

Throughout we use the convention that principal bundles
are right principal bundles.

\begin{definition}\label{def:CartanConnection}
Let \(H \subset G\) be a closed subgroup of a Lie group,
with Lie algebras \(\LieH \subset \LieG\). A \(G/H\)-geometry, or
\emph{Cartan geometry}
modelled on \(G/H\), on a manifold \(M\) is a choice of
$C^\infty$ principal $H$-bundle $E \to M$, and smooth
1-form $\omega \in
\nForms{1}{E} \otimes \LieG$ called the \emph{Cartan connection},
which satisifies all of the following conditions:
\begin{enumerate}
\item
\(
r_h^* \omega = \Ad_h^{-1} \omega
\) for all \(h \in H\).
\item
$\omega_e : T_e E \to \LieG$ is a linear isomorphism at each point
$e \in E$.
\item
For each $A \in \LieG$, define a vector field $\vec{A}$ on $E$ by
the equation $\vec{A} \hook \omega = A$. Then the vector fields
$\vec{A}$ for $A \in \LieH$ generate the $H$-action on $E$.
\end{enumerate}
\end{definition}
Sharpe \cite{Sharpe:2002} gives an introduction to Cartan geometries.
\begin{example}
The principal $H$-bundle $G \to G/H$ is a Cartan geometry, with
Cartan connection $\omega=g^{-1} \, dg$ the left invariant
Maurer--Cartan 1-form on $G$; this geometry is called the
\emph{model Cartan geometry}.
\end{example}
\begin{definition}
An \emph{isomorphism} of $G/H$-geometries $E_0 \to M_0$ and
$E_1 \to M_1$ with Cartan connections $\omega_0$ and $\omega_1$
is an $H$-equivariant diffeomorphism $F : E_0 \to E_1$ 
so that $F^* \omega_1 = \omega_0$.
\end{definition}

\subsubsection{Definition of lift of Cartan geometries}
\begin{definition}
Suppose that $H \subset H' \subset G$ are two 
closed subgroups, and $E \to M'$ is a $G/H'$-geometry.
Let $M=E/H$. Clearly $E \to M$ is a principal $H$-bundle.
We can equip $E$ with the Cartan connection of the
original $G/H'$-geometry, and then clearly $E \to M$ is a $G/H$-geometry.
Moreover $M \to M'$ is a fiber bundle with fiber $H'/H$.
The geometry $E \to M$ is called the \emph{$G/H$-lift} of $E \to M'$ (or simply
the \emph{lift}). Conversely, we will say that a given $G/H$-geometry 
\emph{drops} to a certain
$G/H'$-geometry if it is isomorphic to the lift of that $G/H'$-geometry.
\end{definition}

A Cartan geometry which drops can be completely recovered (up to isomorphism)
from anything it drops to. So dropping encapsulates the same geometry in a lower
dimensional reformulation.

\subsubsection{Definition of generalized flag varieties}

\begin{definition}[Knapp \cite{Knapp:2002}]
A \emph{parabolic subgroup} $P$ of a complex semisimple Lie group $G$
is a subgroup containing a maximal solvable subgroup.
\end{definition}

\begin{remark}
Parabolic subgroups are closed connected complex Lie subgroups.
\end{remark}

\begin{definition}[Landsberg \cite{Landsberg:2005}]
A \emph{generalized flag variety} is a homogeneous
space $G/P$ where $G$ is a complex semisimple Lie group
and $P$ is a parabolic subgroup. Every generalized
flag variety is compact and connected.
\end{definition}

\subsubsection{Semicanonical modules}

\begin{definition}
Suppose that $G/H$ is a complex homogeneous space.
An $H$-submodule $I \subset \left(\LieG/\LieH\right)^*$
is \emph{semicanonical} if there are integers $p \ge 0$ and $q > 0$ so 
that $\left(\det I\right)^{\otimes q} = \left(\det \left(\LieG/\LieH\right)\right)^{\otimes (-p)}$.
An $H$-submodule $I \subset \left(\LieG/\LieH\right)^*$
is \emph{nontrivial} if $I \ne 0$ and $I \ne \left(\LieG/\LieH\right)^*$.
\end{definition}

\subsection{The main theorem}
In various examples, we will prove rigidity of various
Cartan geometries.
Among many other results, we will prove the following theorem:

\begin{theorem}
Suppose that $G$ is a complex Lie group
and $H \subset G$ is a maximal complex subgroup.
Suppose that $I \subset \left(\LieG/\LieH\right)^*$
is a nontrivial semicanonical module.
Suppose that $M$ is a connected smooth complex projective variety
bearing a holomorphic Cartan geometry $E \to M$ modelled on $G/H$.
Let $\mathcal{I} = E \times_H I \subset T^*M$.
If the holomorphic subbundle 
$\mathcal{I}^{\perp} \subset TM$ is not everywhere bracket closed,
then 
\begin{enumerate}
\item
$M=G/H$ and 
\item
the Cartan geometry on $M$ is the model 
holomorphic Cartan geometry on $G/H$ and
\item
$G/H$ is a generalized flag variety.
\end{enumerate}
\end{theorem}

\begin{remark}
We will prove more general theorems below,
and prove a similar theorem for compact K{\"a}hler
manifolds for a different class of Cartan geometries.
We will apply our theorems to prove rigidity
of various types of holomorphic geometric structures. 
\end{remark}

\section{Pfaffian systems}

\begin{definition}
A \emph{Pfaffian system} on a complex manifold $M$ 
is a holomorphic vector subbundle of the holomorphic 
cotangent bundle $T^* M$.
\end{definition}

\begin{remark}
If $\mathcal{I} \subset T^*M$ is a Pfaffian system, 
the reader may feel more comfortable working with
$\mathcal{V}=\mathcal{I}^{\perp}$, which is a holomorphic plane field
(a.k.a. distribution, a.k.a. subbundle of the tangent bundle). 
The convenience of working
with $\mathcal{I}$ rather than $\mathcal{V}$ will become clear, and will
more than overcome the initial discomfort.
\end{remark}

\begin{definition}
A Pfaffian system $\mathcal{I} \subset T^*M$ is \emph{Frobenius} if 
the ideal it generates in the sheaf $\Lm{*}{T^*M}$ of differential
forms is $d$-closed.
\end{definition}

\begin{remark}
Equivalently, $\mathcal{I}$ is Frobenius if 
$\mathcal{V}=\mathcal{I}^{\perp}$ is bracket closed.
Synonyms for \emph{Frobenius} include \emph{integrable, 
completely integrable} and \emph{involutive}.
\end{remark}

\section{Brackets in Cartan geometries}

\begin{lemma}[Sharpe \cite{Sharpe:1997} p. 188, theorem 3.15]\label{lemma:TgtBundle}
If $\pi : E \to M$ is any Cartan geometry,
say with model $G/H$, then the Cartan connection of $E$ maps
\[
\xymatrix{%
0 \ar[r] & \ker \pi'(e) \ar[r] \ar[d] & T_e E \ar[r] \ar[d] & T_m M \ar[r] \ar[d] & 0 \\
0 \ar[r] & \mathfrak{h} \ar[r] & \mathfrak{g} \ar[r] &
\mathfrak{g}/\mathfrak{h} \ar[r] & 0 }
\]
for any points $m \in M$ and $e \in E_m$; thus
\[
TM=E \times_H \left(\mathfrak{g}/\mathfrak{h}\right)
\text{ and }
T^*M=E \times_H \left(\mathfrak{g}/\mathfrak{h}\right)^*.
\]
Under this identification, vector fields on $M$ are identified
with $H$-equivariant functions $E \to \mathfrak{g}/\mathfrak{h}$,
and sections of the cotangent bundle with $H$-equivariant
functions $E \to \left(\mathfrak{g}/\mathfrak{h}\right)^*$.
\end{lemma}

\begin{remark}
If $I \subset \left(\LieG/\LieH\right)^*$ is an $H$-submodule,
then each local holomorphic section of $E \times_H I \subset T^*M$
is identified with an $H$-equivariant holomorphic
map from an open subset of $E$ to $I$.
\end{remark}

\section{Pseudoeffective line bundles}

\begin{definition}
A line bundle $L$ on a K\"ahler manifold is \emph{pseudoeffective}
if $c_1\left(L\right)$ can be represented by a closed positive
$(1,1)$-current. (See Demailly \cite{Demailly:1997a} for more
information.) 
\end{definition}

\begin{remark}
Zero is considered positive in this definition.
\end{remark}

\begin{definition}
If $V$ is a holomorphic vector bundle of rank $N$ on a complex 
manifold $M$, let $\det V = \Lm{N}{V}$.
\end{definition}

\begin{lemma}\label{lemma:Demailly}
Let $M$ be a closed K\"ahler manifold and 
$\mathcal{I} \subset T^*M$  a holomorphic Pfaffian system
which is not Frobenius. Suppose 
that the line bundle $\det \mathcal{I}$ on $M$ 
is pseudoeffective. Then $\mathcal{I}$ is Frobenius.
\end{lemma}
\begin{proof}
Suppose that $\mathcal{I}$ has rank $q$.
Define a line-bundle-valued differential form $\vartheta
\in \nForms{q}{M} \otimes \det\left(TM/\mathcal{I}^{\perp}\right)$ by
\[
\vartheta\left(v_1,v_2,\dots,v_q\right)
=
\left(v_1+I^{\perp}\right) 
\wedge 
\left(v_2+I^{\perp}\right) 
\wedge 
\dots
\wedge 
\left(v_q + I^{\perp}\right).
\]
By Demailly \cite{Demailly:2002} p. 1, Main Theorem applied to
$\vartheta$, if $\det \mathcal{I}$ is pseudoeffective,  then $I$ is Frobenius.
\end{proof}

\begin{definition}
We write the canonical bundle of a complex manifold $M$ as $\kappa_M$.
\end{definition}

\begin{definition}
A holomorphic vector bundle $\mathcal{I}$ on a complex manifold $M$ 
is \emph{semicanonical} if there are integers $p \ge 0, q > 0$
so that $\left(\det \mathcal{I}\right)^{\otimes q} \otimes \kappa_M^{-\otimes p}$
is pseudoeffective.
\end{definition}

\begin{proposition}\label{proposition:ContracanonicalNotPseudoEff}
Suppose that 
\begin{enumerate}
\item $M$ is a compact K{\"ahler} manifold and
\item $\mathcal{I} \subset T^*M$ is a Pfaffian system and
\item $\mathcal{I}$ is not Frobenius and
\item $\mathcal{I}$ is semicanonical.
\end{enumerate}
Then the canonical bundle of $M$ is not pseudoeffective.
\end{proposition}
\begin{proof}
By lemma~\vref{lemma:Demailly}, $\det \mathcal{I}$
is not pseudoeffective. If $\kappa_M$ is pseudoeffective,
then so is $\kappa_M^{\otimes p}$ for any integer $p \ge 0$.
Therefore $\left(\det \mathcal{I}\right)^{\otimes q}$ is pseudoeffective
for some integer $q>0$, and so $\det \mathcal{I}$ is also pseudoeffective.
\end{proof}

\begin{lemma}
Suppose that $\mathcal{I} \subset T^*M$ is a holomorphic
contact structure. Then $\mathcal{I}$ is semicanonical.
\end{lemma}
\begin{proof}
Pick a local section $\vartheta$ of $I$
for which $\vartheta \wedge d \vartheta^n \ne 0$.
But $\vartheta \wedge d \vartheta^n$ is a holomorphic
volume form. The map $\phi : \mathcal{I} \to \kappa_M$,
given on each local section by
$\vartheta \mapsto \vartheta \wedge d\vartheta^n$
depends only on the value of $\vartheta$ pointwise,
and scales like $\phi(f\vartheta)=f^{n+1} \phi(\vartheta)$,
so $\mathcal{I}^{\otimes(n+1)} = \kappa_M$.
\end{proof}

\begin{example}
Suppose that $M=M_1 \times M_2$ is a product.
If $\mathcal{I}_1 \subset T^*M_1$ and $\mathcal{I}_2 \subset T^*M_2$ 
are semicanonical, then so is $\mathcal{I}_1 \oplus \mathcal{I}_2 \subset T^*M$.
\end{example}

\begin{example}\label{example:CartanTwoPlaneField}
Let $V$ be a rank 2 holomorphic subbundle $V \subset TM$
on a complex manifold $M$ with $\dim_{\C{}} M = 5$.
Say that a pair $X, Y$ of local holomorphic sections of $V$ 
is \emph{nondegenerate} if the vector fields
\[
X, Y, [X,Y], \left[X,\left[X,Y\right]\right], \left[Y,\left[X,Y\right]\right]
\] 
are linearly independent at every point where $X$
and $Y$ are defined.
Then $V$ is \emph{Cartan} or \emph{nondegenerate}
if near each point of $M$ there is a nondegenerate pair of 
local holomorphic sections.

 Given any nondegenerate pair $X$ and $Y$,
let
\[
\xi(X,Y) = 
X \wedge Y \wedge [X,Y] \wedge \left[X,\left[X,Y\right]\right] \wedge \left[Y,\left[X,Y\right]\right].
\]
So $\xi$ takes a pair of sections to a section of 
the anticanonical bundle $\kappa_M^*$.
If $X'$ and $Y'$ are any two local sections of $V$,
we can write
\begin{align*}
X' &= a \, X + b \, Y \\
Y' &= c \, X + d \, Y
\end{align*}
for some holomorphic matrix valued function
\[
g = 
\begin{pmatrix}
a & b \\
c & d
\end{pmatrix}
\]
on the overlap where $X', Y', X$ and $Y$ are defined.
Check that
\[
\xi(X',Y')=\det(g)^5 \xi(X,Y).
\]
Therefore if $V$ is nondegenerate,
then $\det V = \Lm{2}{V}$ and $\Lm{2}{V}^{\otimes 5}=\kappa_M^*$. 
Let $\mathcal{I}  = V^{\perp}$.
Clearly $\det \mathcal{I} = \kappa_M \otimes \det V$,
so $\left(\det \mathcal{I}\right)^{\otimes 5} = \kappa_M^{4}$,
so $\mathcal{I}$ is semicanonical.
By Demailly's theorem, since $V$ is not bracket closed 
(i.e. $\mathcal{I}$ is not Frobenius), $\kappa_M$ is pseudoeffective.
\end{example}

\begin{example}\label{example:FirstOrderNondegenerate}
A holomorphic $k$-plane field $V \subset TM$ on a complex manifold
$M$ of dimension $n=\dim_{\C{}} M = k(k+1)/2$ is \emph{first order nondegenerate}
if near each point of $M$ there are local holomorphic sections 
\[
X_1, X_2, \dots, X_k
\]
of $V$ so that 
\[
X_1, X_2, \dots, X_k, 
\left[X_1, X_2\right], \left[X_1,X_3\right], \dots, \left[X_{k-1}, X_k\right]
\]
are linearly independent. Clearly any first order nondegenerate
holomorphic $k$-plane field $V$ has associated Pfaffian
system $\mathcal{I}=V^{\perp}$ semicanonical,
with the same argument as for contact structures.
\end{example}

\section{Pseudoeffectivity and Pfaffian systems in Cartan geometries}

\begin{remark}
Clearly if $I$ is a semicanonical submodule then  $\mathcal{I} = G \times_H I$ 
is a semicanonical vector bundle.
\end{remark}

\begin{example}\label{example:ContactSemicanonicalModule}
If $G/H$ has a $G$-invariant contact structure, say $G \times_H I \subset T^*(G/H)$,
then $I \subset \left(\LieG/\LieH\right)^*$ is semicanonical.
\end{example}

\begin{proposition}
Suppose that 
\begin{enumerate}
\item $M$ is a compact K{\"a}hler manifold and
\item $M$ bears a holomorphic
Cartan geometry modelled on a complex homogeneous
space $G/H$ and
\item $I$ is an $H$-submodule $I \subset \left(\LieG/\LieH\right)^*$ and
\item $I$ is semicanonical and
\item $E \times_H I \subset T^*M$ is not Frobenius.
\end{enumerate}
Then $M$ does not have pseudoeffective canonical bundle.
\end{proposition}
\begin{proof}
Clear from proposition~\vref{proposition:ContracanonicalNotPseudoEff}.
\end{proof}

\begin{example}\label{example:GeneralizedFlagVariety}
Suppose that $G/P$ is a generalized flag variety.
Every $P$-submodule $I \subset \left(\LieG/\LieP\right)^*$
is a direct sum of root spaces, since $P$ contains
the Cartan subgroup of $G$. 
Associate to $I$ the set $S=S_I$ of roots whose root space lies in $I$.
Then $S$ is a set of noncompact positive roots.
For example, $\left(\LieG/\LieP\right)^*$ is the
direct sum of all of the root spaces of all of 
the noncompact positive roots.
Pick any root $\alpha \in S$. Pick $\beta$ to be 
either a compact root or a noncompact
positive root. Then either $\alpha+\beta \in S$ or
$\alpha+\beta$ is not a root.
Conversely, if $S$ is any set of roots with 
this property, let $I=I_S$ be the sum
of the root spaces of the roots that lie in $S$.
Then $I \subset \left(\LieG/\LieP\right)^*$ is a 
$P$-submodule.
So we can draw $I$ by drawing the
root lattice of $G$ and indicating somehow
which roots lie in $S$. 

Next we need to test when $I$ is semicanonical.
Let $W_{G/P}$ be the subgroup of the Weyl group of $G$
preserving the noncompact positive roots of $G/P$.
Baston and Eastwood \cite{Baston/Eastwood:1989}
prove that we can identify the weights of $P$ with the 
$W_{G/P}$-invariant weights of $G$.
The weight of $\det I$ is
\[
\sum_{\alpha \in S} \alpha.
\]
In particular, the weight $\omega$ of $\det\left(\LieG/\LieP\right)^*$ 
as a weight of $P$ is the
sum of the noncompact positive roots, say
\[
\omega=\sum_{\alpha \in \Delta^{\text{noncompact}}_+} \alpha
\]
We can then see that $I$ is semicanonical if and only if
\[
\sum_{\alpha \in S} \alpha 
= \frac{p}{q} \omega
\]
for some rational number $0 \le \frac{p}{q} \le 1$,
and nontrivial if and only if $0 < \frac{p}{q} < 1$.
\end{example}

\begin{example}\label{example:GTwoModPOne}
Figure~\vref{fig:GTwoRootLatticeOne}
\begin{figure}
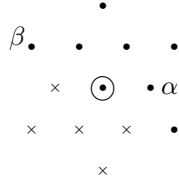

\[
\xy +(2,2)="o",0*\xybox{%
%%0;<3pc,1.5mm>:<0.72pc,1.65pc>::,{"o"
0;<1.5pc,0mm>:<-2.25pc,1.3pc>::,{"o"
\croplattice{-3}3{-2}2{-2.6}{2.6}{-3}3}
,"o"+(1,0)="a"*{\pin}*+!D{}
,"o"+(1.4,0)="aL"*{\alpha}*+!D{}
,"o"+(0,0)="b"*{\bigcirc}*+!L{}
,"o"+(0,0)="c"*{\pin}*+!L{}
,"o"+(0,1)="d"*{\pin}*+!D{}
,"o"+(0,1.2)="dL"*{\beta}*+!D{}
,"o"+(1,1)="B"*{\pin}*+!L{}
,"o"+(2,1)="C"*{\pin}*+!L{}
,"o"+(3,1)="D"*{\pin}*+!L{}
,"o"+(3,2)="E"*{\pin}*+!L{}
,"o"+(-1,0)="na"*{\pout}*+!D{}% \lambda_1}
,"o"+(0,-1)="nd"*{\pin}*+!D{}
,"o"+(-1,-1)="nB"*{\pout}*+!L{}
,"o"+(-2,-1)="nC"*{\pout}*+!L{}
,"o"+(-3,-1)="nD"*{\pout}*+!L{}
,"o"+(-3,-2)="nE"*{\pout}*+!L{}
}
\endxy
\]
\caption{The parabolic subgroup $P_1 \subset G_2$}%
\label{fig:GTwoRootLatticeOne}
\end{figure}
shows the roots of $G_2$. The dots are the roots
whose root spaces lie in $\LieP$, and the crosses
are the other roots. The circled dot
is the origin, representing the Cartan subgroup.
The compact roots lie on the line through $\beta$ and
the origin. The positive noncompact roots lie in the upper half plane
above this line.
The sum of the noncompact positive roots is 
\[
\omega = 10 \, \alpha + 5 \, \beta.
\]

Figure~\vref{fig:GTwoRootLatticeWithRankThreePfaffianSystem} shows
stars ($\star$) on
the noncompact positive roots
\[
3 \, \alpha+2 \, \beta, 2 \, \alpha + \beta, 3 \, \alpha+ \beta.
\]
Let $S$ be the set of these roots, and $I=I_S$.
Then $I$ is a 3-dimensional submodule $I \subset \left(\LieG/\LieP\right)^*$.
\begin{figure}
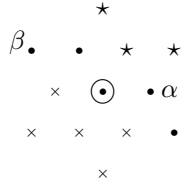

\[
\xy +(2,2)="o",0*\xybox{%
0;<1.5pc,0mm>:<-2.25pc,1.3pc>::,{"o"
\croplattice{-3}3{-2}2{-2.6}{2.6}{-3}3}
,"o"+(1,0)="a"*{\pin}*+!D{}
,"o"+(1.4,0)="aL"*{\alpha}*+!D{}
,"o"+(0,0)="b"*{\bigcirc}*+!L{}
,"o"+(0,0)="c"*{\pin}*+!L{}
,"o"+(0,1)="d"*{\pin}*+!D{}
,"o"+(0,1.2)="dL"*{\beta}*+!D{}
,"o"+(1,1)="B"*{\pin}*+!L{}
,"o"+(2,1)="C"*{\star}*+!L{}
,"o"+(3,1)="D"*{\star}*+!L{}
,"o"+(3,2)="E"*{\star}*+!L{}
,"o"+(-1,0)="na"*{\pout}*+!D{}% \lambda_1}
,"o"+(0,-1)="nd"*{\pin}*+!D{}
,"o"+(-1,-1)="nB"*{\pout}*+!L{}
,"o"+(-2,-1)="nC"*{\pout}*+!L{}
,"o"+(-3,-1)="nD"*{\pout}*+!L{}
,"o"+(-3,-2)="nE"*{\pout}*+!L{}
}
\endxy
\]
\caption{The parabolic subgroup $P_2 \subset G_2$
together with a rank 3 Pfaffian system}%
\label{fig:GTwoRootLatticeWithRankThreePfaffianSystem}
\end{figure}
The weight of $\det I$ is
\begin{align*}
\sum_{\alpha \in S} \alpha
&=
\left(3 \, \alpha+2 \, \beta\right)
+\left(2 \, \alpha + \beta\right) + \left(3 \, \alpha+ \beta\right) 
\\
&= 8 \, \alpha + 4 \, \beta.
\end{align*}
So $I$ is semicanonical.
\end{example}

\begin{example}
Inspect the root lattices of all simple Lie groups $G$ of rank 2, using
the same approach as the previous example.
You see that
for all generalized flag varieties $G/P$ with $G$ of rank 2,
all submodules $I \subset \left(\LieG/\LieP\right)^*$
are semicanonical, except for 
a few counterexamples. These counterexamples
only occur for those $G/P$ where $P=B$
is a Borel subgroup. Specifically $G/P=\SO{5,\C{}}/B$
and $G/P=G_2/B$ have no nontrivial semicanonical modules $I \subset \left(\LieG/\LieP\right)^*$
(i.e. other than $I=0$ and $I=\left(\LieG/\LieP\right)^*$).
On the other hand, $G/P=A_2/B=\SL{3,\C{}}/B$ 
has precisely one nontrivial semicanonical submodule
(as we will see in example~\vref{example:An}),
and various nonsemicanonical submodules.
\end{example}

\begin{example}
More generally, if $B \subset G$ is a
Borel subgroup of a complex semisimple
Lie group $G$, let $V=\LieG_{-\alpha}$
be the root space of any simple root $\alpha$.
Let $I=V^{\perp} \subset \left(\LieG/\LieP\right)^*$,
i.e. $I$ is the sum of the root spaces of all 
positive noncompact roots other than $\alpha$. So the weight of $\det I$
is the sum of all positive noncompact roots other than $\alpha$.
Clearly $I$ is a $B$-submodule of $\left(\LieG/\LieB\right)^*$.
However, $I$ is not semicanonical unless $G=\SL{2,\C{}}$. 
So there are some counterexamples in arbitrary rank.
\end{example}

\begin{proposition}\label{proposition:MaximalGivesSemicanonical}
Suppose that $G$ is a a complex simple Lie group
and $P \subset G$ a maximal parabolic subgroup.
Then every submodule $I \subset \left(\LieG/\LieP\right)^*$ is semicanonical.
\end{proposition}
\begin{proof}
The maximal semisimple subgroup $M \subset P$
from the Langlands decomposition
(see Knapp \cite{Knapp:2002})
has Dynkin diagram given by removing the crossed (i.e.
noncompact) simple roots from the Dynkin diagram of $G/P$. 
Since $P$ is maximal, there is one noncompact simple root, so
the root lattice of $M$ spans a hyperplane
in the root lattice of $G$. All weights of 1-dimensional
$P$-modules lie in the line in the root lattice of $G$
perpendicular to the root lattice of $M$, by invariance
under the Weyl group of $M$.
So if $I$ is a $P$-submodule of $\left(\LieG/\LieP\right)^*$,
then $\det I$ has weight  lying on this line.
The weight $\omega$ 
of $\det \left(\LieG/\LieP\right)^*$ is the sum of the noncompact
positive roots, so is a nonzero vector on this line.
Therefore $\det I$ must have weight a multiple of $\kappa$.
We need to show that this multiple is not negative. This is clear because
the weight is a sum of positive noncompact roots.
\end{proof}

\begin{example}\label{example:An}
The generalized flag variety $G/P$ with $G=A_n$
and Dynkin diagram
\begin{xy}
  0;/r.20pc/: % Set the scale to be quite small
  (10,10)*+{\pout}="2";
  (20,10)*+{\pout}="3";
  (30,10)*+{\pin}="4";
  (40,10)*+{\dots}="5";
  (50,10)*+{\pin}="6";
  "2"; "3" **\dir{-};
  "3"; "4" **\dir{-};
  "4"; "5" **\dir{-};
  "5"; "6" **\dir{-};
\end{xy}
represents the space of pairs $(p,L)$
where $L$ is a projective
line in $\Proj{n}$ and $p \in L$ is a point of that line.
Map $G/P \to \Proj{n}$ by $(p,L) \to p$.
There is an obvious Frobenius Pfaffian system 
$\mathcal{I}_{\text{point}}$
on $G/P$ consisting of the 1-forms vanishing on the fibers of this map.
Similarly there a map $G/P \to \Proj{n*}$, $(p,L) \mapsto L$, and
an obvious Frobenius Pfaffian system 
$\mathcal{I}_{\text{line}}$
on $G/P$ consisting of the 1-forms vanishing on the fibers of this map.
Let $\mathcal{I}_0 = \mathcal{I}_{\text{point}} \cap \mathcal{I}_{\text{line}}$.

As usual, $A_n$ has roots $e_i-e_j \in \R{n+1}$ for $i \ne j$.
A basis of positive simple roots is $\alpha_i = e_i - e_{i+1}$,
$1 \le i \le n$.
The compact roots are $e_i-e_j$ for $i,j \ge 3$ with $i \ne j$.
The noncompact positive roots are $e_1-e_i$ for $i>1$ and $e_2-e_i$
for $i>2$. Write $\alpha \le \beta$ to mean that
$\beta-\alpha$ is a sum of positive noncompact roots and compact roots.

For each positive root $\alpha$, let 
\[
I_{\alpha} = \bigoplus_{\alpha \le \beta} \LieG_{\beta}.
\]
Note that $I_{\alpha} \subset \left(\LieG/\LieP\right)^*$ is a $P$-submodule.
There are precisely 5 distinct $P$-submodules of $\left(\LieG/\LieP\right)^*$:
\begin{enumerate}
\item $0$, 
\item $I_{\text{point}} = I_{\alpha_1}$,  $\dim_{\C{}} I_{\text{point}} = n$,
\item $I_{\text{line}} = I_{\alpha_2}$, $\dim_{\C{}} I_{\text{line}} = 2n-2$,
\item $I_0 = I_{\alpha_1+\alpha_2}$, $\dim_{\C{}} I_0 = n-1$,
and 
\item $\left(\LieG/\LieP\right)^*$.
\end{enumerate}
The associated vector bundles on $G/P$ are the Pfaffian systems 
defined above.

As above let $\omega$ be the sum of the
positive noncompact roots,
\[
\omega = n \, \alpha_1 
+ 2(n-1) \, \alpha_2 
+ 2(n-2) \, \alpha_3 
+
\dots
+ 2 \, \alpha_n,
\]
while the weights for the various submodules are
\begin{align*}
\det I_{\alpha_1}
&: 
n \, \alpha_1
+ (n-1) \, \alpha_2 
+ (n-2) \, \alpha_3 
+ \dots
+ \alpha_n, \\
\det I_{\alpha_2} 
&: 
(n-1) \, \alpha_1
+ 2(n-1) \, \alpha_2 
+ 2(n-2) \, \alpha_3 
+ \dots
+ 2 \, \alpha_n, \\
\det I_0
&:
(n-1) \, \alpha_1
+ (n-1) \, \alpha_2
+ (n-2) \, \alpha_3
+ \dots
+ \alpha_n.
\end{align*}
So $I_0$ is semicanonical
precisely when $n=2$, while $I_{\text{point}}$
and $I_{\text{line}}$ are not semicanonical for any $n$.
\end{example}

\section{Rational curves on smooth complex projective varieties}

\begin{definition}
A complex projective variety is \emph{uniruled}
if every point lies on a rational curve; see \cite{Kollar:1996}.
\end{definition}

\begin{theorem}[Boucksom et. al.%
\cite{Boucksom/Demailly/Paun/Peternell:2004}]\label{theorem:Boucksom}
A smooth complex projective variety is uniruled just when the
variety has nonpseudoeffective canonical bundle.
\end{theorem}

\begin{corollary}\label{corollary:ContainsRationalCurve}
Suppose that $I \subset \left(\LieG/\LieH\right)^*$ is a semicanonical
module. Suppose that $M$ is a smooth complex projective variety 
with a holomorphic Cartan geometry $E \to M$ modelled on $G/H$.
If $E \times_H I$ is not Frobenius then $M$ contains a rational curve.
\end{corollary}

\section{Dropping}

\begin{theorem}%
[Biswas, McKay \cite{Biswas/McKay:2010}]%
\label{thm:CurvesForceDescent}
Suppose that 
\begin{enumerate}
\item 
$G/H$ is a complex homogeneous space,
\item 
$M$ is a connected compact K{\"a}hler manifold and
\item
$M$ bears a holomorphic $G/H$-geometry.
\end{enumerate}
Then the geometry drops to a unique 
\(G/H'\)-geometry on a connected compact
K{\"a}hler manifold $M'$, so that
\begin{enumerate}
\item \(H' \subset G\) is a closed complex subgroup,
\item $H'/H$ is a generalized flag variety,
\item $M \to M'$ is a holomorphic $H/H'$-bundle, and
\item the manifold $M'$ contains no rational curves.
\end{enumerate}
Any other drop $M \to M''$ for which
$M''$ contains no rational curves
factors uniquely through holomorphic drops $M \to M' \to M''$.
\end{theorem}

\begin{theorem}\label{theorem:ContracanonicalDrop}
Suppose that $I \subset \left(\LieG/\LieH\right)^*$ is a semicanonical
module. Suppose that $M$ is a smooth connected 
complex projective variety 
with a holomorphic Cartan geometry $E \to M$ modelled on $G/H$.
Suppose that $E \times_H I$ is not Frobenius.

Then the geometry drops to a unique 
\(G/H'\)-geometry on a connected smooth complex projective
variety $M'$, so that
\begin{enumerate}
\item \(H' \subset G\) is a closed complex subgroup,
\item $H'/H$ is a generalized flag variety,
\item $\dim_{\C{}} H' > \dim_{\C{}} H$, i.e. $\dim_{\C{}} M' < \dim_{\C{}} M$.
\item $M \to M'$ is a holomorphic $H/H'$-bundle, and
\item the manifold $M'$ contains no rational curves.
\end{enumerate}
Any other drop $M \to M''$ for which
$M''$ contains no rational curves
factors uniquely through holomorphic drops $M \to M' \to M''$.

In particular, if there
is no closed proper complex Lie subgroup $H' \subset G$
with $H \subset H'$ and $H'/H$ a rational homogeneous
variety, then $M=G/H$ with its standard flat Cartan geometry,
and $G/H$ is a rational homogeneous variety.
\end{theorem}
\begin{proof}
The manifold $M$ contains a rational curve,
by corollary~\vref{corollary:ContainsRationalCurve}.
By theorem~\vref{thm:CurvesForceDescent},
the geometry drops. If $H$ is not contained
in a closed complex Lie subgroup $H' \subset G$
for which $H'/H$ is a rational homogeneous
variety, then 
the geometry can only drop to a geometry
modelled on $G/G$, a point. The original
geometry on $M$ must be isomorphic to the lift
of $G/G$, i.e. must be isomorphic to $G/H$.
\end{proof}

\begin{definition}
A \emph{parabolic geometry} is a holomorphic Cartan geometry
modelled on a generalized flag variety.
\end{definition}

\begin{remark}
Let's develop a general criterion to ensure that a Pfaffian
system $E \times_P I$ in a parabolic geometry  cannot
be Frobenius. 
There is a well known notion of regularity of
parabolic geometries (see Calderbank and Diemer \cite{Calderbank/Diemer:2001}, 
\v{C}ap \cite{Cap:2006}).
\v{C}ap \cite{Cap:2006} (unnumbered proposition on page 9)
proves that if a parabolic geometry $E \to M$ is regular at a point
of $M$, and if $G \times_P I$ is not Frobenius on $G/P$, 
then $E \times_P I$ is 
also not Frobenius on $M$. We need to see when $G \times_P I$
is Frobenius. It is easy to see that if $I \subset \left(\LieG/\LieP\right)^*$ 
is nontrivial and
semicanonical, then $G \times_P I$ is not Frobenius.
Therefore if $I$ is nontrivial and semicanonical, and
$E \to M$ is regular at a single point of $M$, then $E \times_P I$
is not Frobenius. We will not need to make use of this regularity
criterion in our examples.
\end{remark}

\section{Example: adjoint varieties}

\begin{example}\label{example:adjoint}
Suppose that $G$ is a complex semisimple Lie group.
Pick a highest weight vector $x \in \LieG$, for
some choice of Cartan subalgebra of $G$ and basis of simple roots.
The \emph{adjoint variety} of $G$ is the orbit 
$X=G[x] \subset \Proj{}\LieG$ of the line $[x]$ spanned by $x$ in 
$\LieG$. 
The stabilizer of $[x]$ in $G$ is a parabolic subgroup, 
say $P \subset G$ and $X=G/P$.
For example, if $G$ is simple, the adjoint varieties
have Dynkin diagrams as in figure~\vref{figure:adjointVarieties}.
\begin{figure}
\[
 \begin{array}{llcl}
 \textbf{Group} & \textbf{Variety} & \textbf{dim} & \textbf{Diagram} \\
  A_n & \Proj{}T^*\Proj{n} & 2n-1 &
{\scriptscriptstyle{ \begin{xy}
  0;/r.20pc/: % Set the scale to be quite small
 (0,0)*+{\pout}="1";
  (10,0)*+{\pin}="2";
  (20,0)*+{\dots}="3";
  (30,0)*+{\pin}="4";
  (40,0)*+{\pout}="5";
  "1"; "2" **\dir{-};
  "2"; "3" **\dir{-};
  "3"; "4" **\dir{-};
  "4"; "5" **\dir{-};
  \end{xy} }} \\
  B_n & \Grnull{2}{2n+1} & 4n-5 &
{\scriptscriptstyle{  \begin{xy}
  0;/r.20pc/: % Set the scale to be quite small
  (0,0)*+{\pin}="1";
  (10,0)*+{\pout}="2";
  (20,0)*+{\pin}="3";
  (20,0)*+{\dots}="4";
  (30,0)*+{\pin}="5";
  (40,0)*+{\pin}="6";
  (50,0)*+{\pin}="7";
%  (80,0)*+{\Grnull{2}{2n+1}};
  "1"; "2" **\dir{-};
  "2"; "3" **\dir{-};
  "3"; "4" **\dir{-};
  "4"; "5" **\dir{-};
  "5"; "6" **\dir{-};
  {\ar@{=>}; "6";"7"};
  \end{xy} }} \\
  C_n & \Proj{2n-1} & 2n-1 &
{\scriptscriptstyle{
  \begin{xy}
  0;/r.20pc/: % Set the scale to be quite small
  (0,0)*+{\pout}="1";
  (10,0)*+{\pin}="2";
  (20,0)*+{\dots}="3";
  (30,0)*+{\pin}="4";
  (40,0)*+{\pin}="5";
  (50,0)*+{\pin}="6";
  %(70,0)*+{\Proj{2n-1}};
  "1"; "2" **\dir{-};
  "2"; "3" **\dir{-};
  "3"; "4" **\dir{-};
  "4"; "5" **\dir{-};
  {\ar@{=>}; "6";"5"};
  \end{xy} }} \\
  D_n & \Grnull{2}{2n} & 4n-7 & 
{\scriptscriptstyle{
\begin{xy}
 0;/r.20pc/: % Set the scale to be quite small
(-20,0)*+{\pin}="1"; (-10,0)*+{\pout}="2"; (0,0)*+{\dots}="3";
(10,0)*+{\pin}="4"; (20,0)*+{\pin}="5"; (25,9)*+{\pin}="6";
(25,-9)*+{\pin}="7"; 
%(50,0)*+{\Grnull{2}{2n}}; 
"1"; "2" **\dir{-}; "2";
"3" **\dir{-}; "3"; "4" **\dir{-}; "4"; "5" **\dir{-}; "5"; "6"
**\dir{-}; "5"; "7" **\dir{-};
\end{xy}}} \\
E_6 & X^{\ad}_{E_6} & 21 & {\scriptscriptstyle{
\begin{xy}
 0;/r.20pc/: % Set the scale to be quite small
  (0,0)*+{\pin}="1";
  (10,0)*+{\pin}="2";
  (20,0)*+{\pin}="3";
  (20,10)*+{\pout}="4";
  (30,0)*+{\pin}="5";
  (40,0)*+{\pin}="6";
%  (70,0)*+{\Grw};
  "1"; "2" **\dir{-};
  "2"; "3" **\dir{-};
  "3"; "4" **\dir{-};
  "3"; "5" **\dir{-};
  "5"; "6" **\dir{-};
\end{xy}
}} \\
E_7 & X^{\ad}_{E_7} & 33 & {\scriptscriptstyle{
\begin{xy}
 0;/r.20pc/: % Set the scale to be quite small
  (0,0)*+{\pout}="1";
  (10,0)*+{\pin}="2";
  (20,0)*+{\pin}="3";
  (20,10)*+{\pin}="4";
  (30,0)*+{\pin}="5";
  (40,0)*+{\pin}="6";
  (50,0)*+{\pin}="7";
%  (70,0)*+{\Grw};
  "1"; "2" **\dir{-};
  "2"; "3" **\dir{-};
  "3"; "4" **\dir{-};
  "3"; "5" **\dir{-};
  "5"; "6" **\dir{-};
  "6"; "7" **\dir{-};
\end{xy}
}} \\
E_8 & X^{\ad}_{E_8} & 57 & {\scriptscriptstyle{
\begin{xy}
 0;/r.20pc/: % Set the scale to be quite small
  (0,0)*+{\pin}="1";
  (10,0)*+{\pin}="2";
  (20,0)*+{\pin}="3";
  (20,10)*+{\pin}="4";
  (30,0)*+{\pin}="5";
  (40,0)*+{\pin}="6";
  (50,0)*+{\pin}="7";
  (60,0)*+{\pout}="8";
%  (70,0)*+{\Grw};
  "1"; "2" **\dir{-};
  "2"; "3" **\dir{-};
  "3"; "4" **\dir{-};
  "3"; "5" **\dir{-};
  "5"; "6" **\dir{-};
  "6"; "7" **\dir{-};
  "7"; "8" **\dir{-};
\end{xy}
}} \\
F_4 & X^{\ad}_{F_4} & 15 & {\scriptscriptstyle{
\xy
 0;/r.10pc/: % Set the scale to be quite small
  (10,0)*+{\pout}="0";
  (30,0)*+{\pin}="1";
  (40,0)*+{>}="2";
  (50,0)*+{\pin}="3";
  (70,0)*+{\pin}="4";
  "0"; "1" **\dir{-};
  "1"; "3" **\dir2{-};
  "3"; "4" **\dir{-};
\endxy
}} \\
G_2 & \Grnull{2}{\operatorname{Im} \mathbb{O}} & 5 & {\scriptscriptstyle{
\xy
 0;/r.10pc/: % Set the scale to be quite small
  (30,0)*+{\pin}="1";
  (40,0)*+{<}="2";
  (50,0)*+{\pout}="3";
  "1"; "3" **\dir3{-};
\endxy
}}
\end{array}
\]
\caption{The adjoint varieties of the complex simple Lie groups. The adjoint variety of $C_n$ 
is $\Proj{2n-1}$ under the Veronese embedding, i.e. the set 
of rank 1 quadratic forms up to rescaling. A 2-plane in the
octave numbers $\mathbb{O}$ is \emph{null} if the multiplication is zero on it.
The adjoint variety of $G_2$ is the set 
of null 2-planes in the imaginary complexified octave numbers. We don't know a
geometric description of the adjoint varieties
of $E_6, E_7, E_8$ and $F_4$.}%
\label{figure:adjointVarieties}
\end{figure}

If $G$ is simple, then its adjoint variety is a holomorphic
contact manifold, and every homogeneous compact complex contact
manifold occurs as an adjoint variety; see Landsberg \cite{Landsberg:2005}.
There is precisely one holomorphic contact structure
on any adjoint variety.

If $G$ is not simple, then up to a finite covering $G$ is a product 
of simple factors 
\[
G = G_1 \times G_2 \times \dots \times G_s,
\]
and correspondingly
\[
P = P_1 \times P_2 \times \dots \times P_s,
\]
where $P_j = P \cap G_j$. For each $G_j$, we can then
consider the one dimensional $P_j$-submodule $I_j \subset \left(\LieG_j/\LieP\right)^*$
which arises from the holomorphic contact
structure on $X_j=G_j/P_j$. We can then let $I=\bigoplus_j I_j$,
and again $I$ is semicanonical on $X=G/P$, though not a contact structure.
\end{example}

\begin{theorem}
Suppose that $G$ is a complex simple Lie group
and that $G/P$ is an adjoint variety with holomorphic
contact structure $G \times_P I$.
Suppose that $E \to M$ is holomorphic parabolic geometry 
modelled on $G/P$,
on a smooth connected complex projective
variety $M$. Let $\mathcal{I} = E  \times_P I \subset T^*M$.  Either
\begin{enumerate}
\item $M$ is foliated by smooth hypersurfaces on which $\mathcal{I}=0$ or
\item $M=G/P$ with its usual adjoint variety geometry or
\item $G=A_n$, and the geometry on $M$ drops to a
holomorphic projective connection
on a smooth connected complex projective variety.
\end{enumerate}
\end{theorem}
\begin{proof}
Either $\mathcal{I}$ is Frobenius, or the parabolic geometry drops by
theorem~\vref{theorem:ContracanonicalDrop}.

The adjoint variety $X=G/P$ of $G=A_n$ is the variety of pairs of
a hyperplane in $\Proj{n}$ and a point on that hyperplane.
There are only two parabolic subgroups of $A_n$ containing
$P$: forget the point or the hyperplane, i.e. $G/P'$ is either 
projective space or the dual projective space. Projective
space and its dual are isomorphic, so the same parabolic geometries
are modelled on either one. 
Suppose that $M$ is a smooth complex projective variety with a holomorphic
parabolic geometry modelled on the adjoint variety of $A_n$. 
Then $M$ drops to a smooth complex projective variety 
with holomorphic projective connection.

Consider the adjoint variety $X=G/P$ of any other simple complex Lie group $G$ (i.e.
$G=B_n, C_n, D_n, E_6, E_7, E_8, F_4$ or $G_2$). Then
$P \subset G$ is a maximal parabolic subgroup. So there is 
only one regular parabolic geometry on any smooth complex projective variety modelled
on that $G/P$: the model $G/P$ with its standard flat $G/P$-geometry.
\end{proof}

\begin{remark}
We will reconsider the $A_n$-adjoint geometries in 
section~\vref{section:DoubleLegendreFoliations}.
\end{remark}

\section{Example: Cartan's theory of 2-plane fields on 5-manifolds}

In example~\vref{example:GTwoModPOne}, we saw that
$G_2/P_1$ bears a holomorphic rank 3 Pfaffian system.
We can see from the root lattice in figure~\vref{fig:GTwoRootLatticeOne}
that $\dim_{\C{}} G_2/P_1 = 5$ (i.e. 5 crosses representing
the 5 dimensions of $\LieG_2/\LieP_1$). We can also
see that the rank 3 Pfaffian system is not Frobenius,
because there is a pair of noncompact positive roots not among those
3 which add up to a root among those 3. 
The dual plane field is associated to the $P_1$-module $V=I^{\perp}$,
i.e. the sum of root  spaces of the two roots $-\alpha,-\alpha-\beta$.
We can even see that the
2-plane field is Cartan, in the sense of 
example~\vref{example:CartanTwoPlaneField},
by looking at the brackets of vector fields in $\LieG_2/\LieP_1$,
i.e. looking at sums of the roots $-\alpha,-\alpha-\beta$.
(We leave this claim to the reader to prove, since it is
not essential to our arguments.)

\begin{theorem}[Cartan \cite{Cartan:30,Sternberg:1983,Tanaka:1979}]
If $\mathcal{V}$ is a Cartan 2-plane field on a
5-dimensional complex manifold $M$, then
then there is a holomorphic
parabolic geometry $E \to M$ modelled on $G_2/P_1$,
so that $\mathcal{V} = E \times_P V \subset TM$, where 
$V \subset \LieG_2/\LieP_1$ is the $P_1$-submodule 
constructed in example~\vref{example:GTwoModPOne}.
\end{theorem}

\begin{theorem}
The only holomorphic Cartan 2-plane 
field on any smooth connected complex projective variety
is the standard one on $G_2/P_1$ described
in example~\vref{example:GTwoModPOne}.
\end{theorem}
\begin{proof}
Suppose that $M$ is a smooth connected complex
projective variety of complex dimension 5,
bearing a holomorphic Cartan 2-plane field.
From example~\vref{example:CartanTwoPlaneField}, we have seen
that a smooth complex projective variety with a Cartan 2-plane
field must have nonpseudoeffective canonical bundle.
By theorem~\vref{theorem:Boucksom}, the 
variety must then be uniruled.

By Cartan's theorem, we can assume that the
Cartan 2-plane field is $E \times_P V \subset TM$.
Let $I = V^{\perp}$. Since $I$ is semicanonical, and the Pfaffian system 
$E \times_P I$ is not Frobenius, again we see that 
the variety $M$ must be uniruled.
By theorem~\vref{thm:CurvesForceDescent},
the parabolic geometry must drop to a parabolic
geometry with a lower dimensional model.
The parabolic geometry can only drop 
to a parabolic geometry modelled on a point, since $P$ is maximal,
so drops just when the parabolic geometry is isomorphic to the model.
\end{proof}

\section{Example: 3-plane fields on 6-manifolds}

\begin{definition}
A rank 3 Pfaffian system $\mathcal{I} \subset T^*M$ 
on a complex manifold $M$ of complex dimension 6
is \emph{nondegenerate} if near each point of $M$
there are 3 sections of $\mathcal{I}$ with linearly
independent exterior derivatives.
\end{definition}

\begin{example}\label{example:BryantExample}
Let $G=B_3=\PO{7,\C{}}$ and $G/P$
be the space of null 3-planes in $\C{6}$
for some nondegenerate complex inner product.
The Dynkin diagram of $G/P$ is
  \begin{xy}
  0;/r.10pc/: % Set the scale to be quite small
  (0,0)*+{\pin}="1";
  (22,0)*+{\pin}="2";
  (44,0)*+{\pout}="3";
  (33,0)*+{>}="edge";
%  (33,7)*+{\text{Nondegenerate 3-plane field on 6-fold}};
  "1"; "2" **\dir{-};
  {\ar@{=}; "2";"3"};
\end{xy}.
Write the simple roots of $G$ as $\alpha_1=e_1-e_2, \alpha_2=e_2-e_3, \alpha_3=e_2+e_3$
in terms of the standard basis $e_1, e_2, e_3 \in \R{3}$.
The root $\alpha_3$ will be  the noncompact positive simple root.
The noncompact positive roots are
\[
\alpha_3,
\alpha_2+\alpha_3,
\alpha_1+\alpha_2+\alpha_3,
\alpha_2 + 2 \, \alpha_3, 
\alpha_1 + \alpha_2 + 2 \, \alpha_3,
\alpha_1 + 2 \, \alpha_2 + 2 \, \alpha_3.
\]
Let $S$ be the set of roots 
\[
\alpha_2 + 2 \, \alpha_3, 
\alpha_1 + \alpha_2 + 2 \, \alpha_3,
\alpha_1 + 2 \, \alpha_2 + 2 \, \alpha_3
\]
(i.e. the roots with $2 \, \alpha_3$ in them).
Let $I=I_S \subset \left(\LieG/\LieP\right)^*$
be the sum of the root spaces of roots in $S$;
$I$ has dimension 3. We can
see that $G \times_P I$ is nondegenerate
in Bryant's sense, since we can write the 3 roots
in $S$ each as a sum of distinct pairs of roots
not in $S$. (We again leave the reader to figure out the yoga
relating exterior derivatives to root sums, since we
won't use this fact.)
\end{example}

\begin{theorem}%
[Bryant \cite{Bryant:2005}]%
\label{theorem:Bryant}
If $\mathcal{I} \subset T^*M$ is a nondegenerate
rank 3 Pfaffian system on a smooth complex
projective variety $M$ with $\dim_{\C{}} M = 6$,
then there is a parabolic geometry $E \to M$
so that $\mathcal{I} = E \times_P I$.
\end{theorem}

\begin{theorem}\label{theorem:BryantDrop}
Suppose that $M$ is a smooth connected complex projective
variety of complex dimension 6, bearing a nondegenerate
rank 3 Pfaffian system. Then $M=B_3/P$ is the
model defined in example~\vref{example:BryantExample}.
\end{theorem}
\begin{proof}
Suppose that $\mathcal{I} \subset T^*M$ 
is a nondegenerate rank 3 Pfaffian
system on a 6-dimensional 
connected smooth complex projective variety $M$. 
By Bryant's theorem, we can
assume that $\mathcal{I} = E \times_P I$, for
some parabolic geometry $E \to M$.
Apply theorem~\vref{theorem:ContracanonicalDrop} to prove that 
the geometry on $M$ drops. The group
$P \subset B_3$ is a maximal parabolic subgroup.
Therefore $M$ must drop
to a point, i.e. must be isomorphic to $B_3/P$.
\end{proof}

\section{Example: quaternionic contact structures}%
\label{section:QuaternionicContact}

\begin{definition}
Suppose that $M$ is a complex manifold,
$\dim_{\C{}} M = 7$ and that $\mathcal{I} \subset TM$
is a holomorphic Pfaffian system of rank 3.
For any two local sections $\vartheta_0, \vartheta_1$
of $\mathcal{I}$, let
\[
q\left(\vartheta_0,\vartheta_1\right) = 
\left.
d\vartheta_0 \wedge d\vartheta_1
\right|_{\mathcal{I}^{\perp}}.
\]
It is easy to check that $q$ is a global holomorphic section
of 
\[
\Sym{2}{\mathcal{I}}^* \otimes \Lm{4}{\mathcal{I}^{\perp}}.
\]
Say that $\mathcal{I}$ is \emph{nondegenerate} if
$\vartheta_0 \hook q=0$ precisely when $\vartheta_0=0$.
A \emph{quaternionic contact structure} is a nondegenerate
holomorphic Pfaffian system of rank 3 on a 7-manifold.
\end{definition}

\begin{remark}
Quaternionic contact structures are very clearly 
explained by Montgomery \cite{Montgomery:2002}.
For discussion of real forms of quaternionic contact
structures, see \cite{Biquard:2000,Biquard:2001,Fox:2005}.
\end{remark}

\begin{example}\label{example:StandardQCS}
Let $X=C_3/P=\Sp{6,\C{}}/P$ the space of subLagrangian
2-planes in $\C{6}$, where $P$ is the stabilizer of a subLagrangian
2-plane.  The Dynkin diagram of $X$ is 
  \begin{xy}
  0;/r.10pc/: % Set the scale to be quite small
  (0,0)*+{\pin}="1";
  (22,0)*+{\pout}="2";
  (44,0)*+{\pin}="3";
  (33,0)*+{<}="edge";
  "1"; "2" **\dir{-};
  {\ar@{=}; "2";"3"};
\end{xy}
There is a $C_3$-invariant quaternionic contact 
structure on $X$ defined as follows. 

We can write the roots of $C_3$ as vertices
and the middles of edges of an octahedron,
say as $\pm e_i \pm e_j$ for $1 \le i,j \le 3$.
The positive simple roots are
\[
\alpha_1 = e_1 - e_2, \alpha_2 = e_2 - e_3, \alpha_3 = 2 e_3.
\]
The positive noncompact roots of $X$ are
\[
\alpha_2,  \, \alpha_1+\alpha_2, \, \alpha_2+\alpha_3, \, \alpha_1 + \alpha_2 + \alpha_3, \, 
2 \alpha_2 + \alpha_3, \,  \alpha_1  + 2 \, \alpha_2 + \alpha_3, \,
2 \, \alpha_1  + 2 \, \alpha_2 + \alpha_3.
\]
Consider the 3 roots
\[
2 \alpha_2 + \alpha_3 , \, \alpha_1  + 2 \, \alpha_2 + \alpha_3, \,
2 \, \alpha_1  + 2 \, \alpha_2 + \alpha_3,
\]
i.e. those with $2 \alpha_2$ in them.
Let $I \subset \LieP$ be the sum of the root spaces of those
3 roots. Consider
the Pfaffian system $\mathcal{I} = C_3 \times_P I \subset T^*\left(C_3/P\right)$.
One can directly calculate using the structure equations
of $C_3$ that $\mathcal{I}$ is a quaternionic contact structure.
(Once more we leave this local calculation to the reader, since
we don't need this result.)
\end{example}

\begin{theorem}%
[Montgometry \cite{Montgomery:2002}]
If $\mathcal{I}$ is a quaternionic contact structure
on a complex manifold $M$, then there is
a holomorphic parabolic geometry $E \to M$ modelled
on $X=C_3/P$, so that $E \times_P I = \mathcal{I}$,
for $I = \left(\LieG/\LieP\right)^*$ the semicanonical $P$-submodule
defined in example~\vref{example:StandardQCS}.
\end{theorem}

We now prove theorem~\vref{theorem:RigidityForQuaternionicContact}.
%\begin{theorem}
%Suppose that $M$ is a smooth connected complex projective variety
%of complex dimension 7, bearing a quaternionic contact
%structure. Then $M$ is the model $C_3/P$ with its standard
%flat quaternionic contact structure.
%\end{theorem}
\begin{proof}
Suppose that $\mathcal{I}$ is a quaternionic contact
structure on $M$. By Montgomery's theorem, we can assume that
$\mathcal{I} = E \times_P I$ for some holomorphic
parabolic geometry $E \to M$.
Apply theorem~\vref{theorem:ContracanonicalDrop} to prove that 
the parabolic geometry on $M$ drops. Since $P \subset C_3$ is a 
maximal parabolic subgroup,
$M$ must drop
to a point, i.e. must be isomorphic to $X$.
\end{proof}

\section{Example: \v{C}ap--Neusser Pfaffian systems}

\begin{example}\label{example:CapNeusserTwo}
For $n \ge 3$, let $G=B_n=\PO{2n+1,\C{}}$, Let $X=G/P$
be the set of all $n$-dimensional null subspaces
of the standard complex linear inner product
on $\C{2n+1}$. The Dynkin diagram of $X$
is
\[
{\scriptstyle{
  \begin{xy}
  0;/r.20pc/: % Set the scale to be quite small
  (0,0)*+{\pin}="1";
  (10,0)*+{\pin}="2";
  (20,0)*+{\dots}="3";
  (30,0)*+{\pin}="4";
  (40,0)*+{\pin}="5";
  (50,0)*+{\pout}="6";
%  %(70,0)*+{Q^{2n-1}};
  "1"; "2" **\dir{-};
  "2"; "3" **\dir{-};
  "3"; "4" **\dir{-};
  "4"; "5" **\dir{-};
  {\ar@{=>}; "5";"6"};
  \end{xy}
}}
\] 
If we order the positive roots according
to the coefficient of $\alpha_n$,
there are precisely $n$ positive noncompact
roots $\alpha$ with coefficient 1
and precisely $n(n-1)/2$ positive noncompact
roots $\alpha$ with coefficient 2.
Let $S$ be the set of noncompact positiive 
roots of coefficient 2. Let $I=I_S$ be the
sum of the root spaces of these roots, 
so $I \subset \left(\LieG/\LieP\right)^*$. 
Then $I$ turns out to be a first order
nondegenerate Pfaffian system in the
sense of example~\vref{example:FirstOrderNondegenerate}.
(Again we leave this statement for the reader
to prove.)
\end{example}

\begin{theorem}%
[\v{C}ap and Neusser \cite{Cap/Neusser:2009}]%
\label{theorem:CapNeusser}
Suppose that $n \ge 3$.
Suppose that $M$ is a complex manifold with
$\dim_{\C{}} M = n(n+1)/2$.
Suppose that $\mathcal{I} \subset T^*M$
is a holomorphic first order nondegenerate Pfaffian system.
Then there is a holomorphic parabolic geometry
$E \to M$ so that $E \times_P I=\mathcal{I}$,
where $I$ is the semicanonical $P$-module
defined in example~\vref{example:CapNeusserTwo}.
\end{theorem}

\begin{theorem}\label{theorem:CapNeusserDrop}
Suppose that $M$ is a smooth connected
complex projective variety with
$\dim_{\C{}} M \ge 6$ bearing
a holomorphic first order nondegenerate Pfaffian system.
Then $M=B_n/P$ with its standard first
order nondegenerate Pfaffian system as 
defined in example~\vref{example:CapNeusserTwo}.
\end{theorem}
\begin{proof}
Suppose that $\mathcal{I}$ is a first order nondegenerate
Pfaffian system on $M$.
By theorem~\vref{theorem:CapNeusser} of \v{C}ap 
and Neusser, we can assume
that $\mathcal{I} = E \times_P I$ for some
holomorphic parabolic geometry $E \to M$
modelled on $B_n/P$.
Apply theorem~\vref{theorem:ContracanonicalDrop} to prove that 
the geometry on $M$ drops. Since the model $P \subset B_n$
is a maximal parabolic subgroup, $M$ must drop
to a point, i.e. must be isomorphic to $X$.
\end{proof}

\begin{remark}
Theorem~\vref{theorem:BryantDrop} is the special case of
theorem~\vref{theorem:CapNeusserDrop}
for $\dim_{\C{}} M = 6$.
\end{remark}

\section{Example: parabolic geometries modelled on products}

\begin{theorem}
Suppose that $E \to M$ is a holomorphic parabolic
geometry on a smooth complex projective variety $M$,
modelled on a generalized flag variety $G/P$. 
Suppose that $G$ splits into a product of simple complex Lie groups,
\[
G = G_1 \times G_2 \times \dots \times G_s.
\]
Let $P_j = P \cap G_j$ for each $j$. 
Suppose that $P_j \subset G_j$ is maximal for each $j$.
Suppose that each $G_j/P_j$
has a nontrivial semicanonical $P_j$-submodule 
$I_j \subset \left(\LieG_j/\LieP_j\right)^*$.
Let $I=I_1 \oplus I_2 \oplus \dots \oplus I_s$.
Suppose that every $E \times_P I_j$ is not Frobenius.
Then $M=G/P$ with its standard flat parabolic
geometry.
\end{theorem}
\begin{proof}
Theorem~\ref{theorem:ContracanonicalDrop}
ensures that the parabolic geometry drops,
say to a geometry with some model $G/Q$, $P \subset Q \subset G$
on some complex manifold $M'$.
Since each $P_j \subset G_j$ is maximal, 
the group $Q$ must be obtained by setting
$Q_j=P_j$ or $Q_j=G_j$ for each value of $j$,
and then $Q = Q_1 \times Q_2 \times \dots \times Q_s$.
If $P_j \ne Q_j$, then $E \times_Q I_j$ is 
not Frobenius, since its local sections
pull back to local sections of $E \times_Q I_j$.
Therefore $Q=G$,
and therefore $M'$ is a point, and so $M$ must
be isomorphic to the model, i.e. $M=G/P$.
\end{proof}

\section{Example: double Legendre foliations}%
\label{section:DoubleLegendreFoliations}

We arrive at our most complicated example.
We will study parabolic geometries modelled
on the adjoint variety of $A_n$, but we will
not obtain a complete classification.

\begin{example}\label{example:AnAdjointAsDLF}
Consider the adjoint variety of $A_n=\SL{n+1,\C{}}$,
say $X=G/P$, $G=A_n$. 
Let's first find all of the $P$-submodules $I \subset \left(\LieG/\LieP\right)^*$.
Write  the positive roots of $A_n$
as $\alpha_i + \alpha_{i+1} + \dots + \alpha_j$
for $i \le j$. Let $S_1$ be the set of roots
\[
\alpha_1, \alpha_1 + \alpha_2, \dots, 
\alpha_1 + \alpha_2 + \dots + \alpha_n
\]
Let $S_n$ be the set of roots
\[
\alpha_1 + \alpha_3 + \dots + \alpha_n,
\alpha_2 + \alpha_4 + \dots + \alpha_n,
\dots,
\alpha_n.
\]
Let $I_1 \subset \left(\LieG/\LieP\right)^*$ 
be the sum of all root spaces of positive noncompact
roots $\alpha$ for which $\alpha \in S_1$, and similarly
let $I_n \subset \left(\LieG/\LieP\right)^*$ be 
the sum of all root spaces of positive noncompact
roots $\alpha$ for which $\alpha \in S_n$.
It is clear that no two noncompact positive roots can
add up to a root in $S_1$, and similarly for $S_n$.
It turns out to follow that $G \times_P I_1$ and 
$G \times_P I_n$ are Frobenius. 
(We leave the reader to figure out the yoga
relating exterior derivatives to root sums, since
we will only make use of it in examples where
the claims made are an elementary calculation
using Cartan's structure equations.)
By a similar argument, if we let $I_{1n} = I_1 \cap I_n$,
then $G \times_P I_{1n}$ is a contact structure.
Indeed $I_{1n}$ corresponds to the set $S_{1n}=S_1 \cap S_n$,
which is just
\[
\alpha_1 + \alpha_2 + \dots + \alpha_n,
\]
which is a sum of noncompact positive roots
\[
\left(\alpha_1 + \alpha_2 + \dots + \alpha_{j-1}\right)
+
\left(\alpha_j + \alpha_{j+1} + \dots + \alpha_n\right),
\]
corresponding to the exterior derivative of the 
contact form being a sum of the wedge products
of various 2-forms. 
The leaves of $G \times_P I_1$ are the fibers of $G/P \to \Proj{n}$,
while the leaves of $G \times_P I_n$ are the
fibers of $G/P \to \Proj{n*}$.  In particular,
$G/P$ has two foliations (indeed fiber bundle mappings),
with leaves integral manifolds of the contact structure.
\end{example}

\begin{definition}%
[Tabachnikov \cite{Tabachnikov:1993}]%
A \emph{double Legendre foliation} \cite{Tabachnikov:1993} of a complex
manifold $M$ of complex dimension $2n-1$
is a pair $F_0,F_1 \subset TM$ of holomorphic foliations 
so that $F_0 \oplus F_1 \subset TM$ is a holomorphic contact structure,
and the leaves of $F_0$ and of $F_1$ are Legendre submanifolds.
\end{definition}

\begin{example}
The $A_n$-adjoint variety has a holomorphic double Legendre foliation.
\end{example}

\begin{theorem}%
[Tabachnikov \cite{Tabachnikov:1993}]%
\label{theorem:Tabachnikov}
Suppose that $F_0, F_1$ is a holomorphic 
double Legendre foliation of a complex manifold $M$
of complex dimension $2n-1$.
Then there is a  holomorphic parabolic geometry $E \to M$ modelled
on the adjoint variety of $A_n=\SL{n+1,\C{}}$ so that (in the notation of example~\vref{example:AnAdjointAsDLF})
$F_0 = E \times_P I_1^{\perp}$ and $F_1 = E \times_P I_n^{\perp}$.
\end{theorem}

\begin{example}\label{example:AdjointAn}
We will construct a parabolic geometry modelled
on the adjoint variety of $A_n$ by lifting a holomorphic
projective connection. 

Write points of $\C{n+1}$
as columns, spanned by the standard basis $e_0, e_1, \dots, e_n$.
Clearly $\Proj{n}=\PSL{n+1,\C{}}/P_1$
where $P_1 \subset G=\PSL{n+1,\C{}}$ is the subgroup 
of matrices of the form
\[
\begin{bmatrix}
p^0_0 & p^0_j \\
0 & p^i_j
\end{bmatrix},
\]
$1 \le i,j \le n$. 
and we identify an element of $\LieG/\LieP_1$
with a column in $\C{n}$ by writing out the entries 
\[
\begin{pmatrix}
A^1_0 \\
A^2_0 \\
\vdots \\
A^n_0
\end{pmatrix}
\in \C{n}.
\]

Let $P_n \subset G=\PSL{n+1,\C{}}$ be the subgroup 
of matrices of the form
\[
\begin{bmatrix}
p^0_0 & p^0_j & p^0_n \\
p^i_0 & p^i_j & p^i_n \\
0 & 0 & p^n_n
\end{bmatrix},
\]
$1 \le i,j \le n-1$. 
Let $P =P_1 \cap P_n \subset G=\PSL{n+1,\C{}}$ be the subgroup 
of matrices of the form
\[
\begin{bmatrix}
p^0_0 & p^0_j & p^0_n \\
0 & p^i_j & p^i_n \\
0 & 0 & p^n_n
\end{bmatrix},
\]
$1 \le i,j \le n-1$. 

Suppose that $M'$ is a complex manifold of complex
dimension $n$ bearing a holomorphic projective connection,
i.e. a holomorphic parabolic geometry $\pi : E \to M'$
modelled on $\Proj{n}$. We will construct the lift
of $M'$ to a parabolic geometry modelled on the adjoint
variety of $A_n$. Let $M=\Proj{}T^*M'$, with its usual holomorphic
contact structure.  
Write points of $M$ as $m=\left(m',H\right)$
where $m' \in M'$ and $H \subset T_{m'} M'$
is a complex hyperplane. Map $E \to M$ by
\[
e \in E \mapsto m=\left(m',H\right) \in M,
\] 
taking $H$ to be the hyperplane identified
by the Cartan connection $\omega$ with the 
span of $e_1, e_2, \dots, e_{n-1} \in \C{n}=\LieG/\LieP_1$.
Because $\omega$ transforms in the adjoint $P$-representation,
i.e.
\[
r_p^* \omega = \Ad(p)^{-1} \omega,
\]
we can easily check that $M=E/P$. Therefore
$M=\Proj{}T^*M'$ is the lift of $M'$, and $E \to M$
is a parabolic geometry
modelled on the adjoint variety of $A_n$.

Consider on $E$ the following two linear Pfaffian systems:
let $\mathcal{I}_0 \subset T^* E$ be the system
\[
\omega+\LieP_1 = 0,
\]
and let $\mathcal{I}_1 \subset T^*E$ be the system
\[
\omega + \LieP_n = 0.
\]
The fibers of $E \to M$ are Cauchy characteristics
for each of these systems, and both systems are 
$P$-invariant. Therefore these Pfaffian
systems are pulled back from Pfaffian systems,
which we denote by the same names, on $M$; see \cite{BCGGG:1991}.
Clearly on $M$, $\mathcal{I}_0 = E \times_P I_0$
and $\mathcal{I}_1 = E \times_P I_1$.
Let $F_0=\mathcal{I}_0^{\perp}$ and $F_1=\mathcal{I}_1^{\perp}$.

This holomorphic vector subbundle $F_1 \subset TM$ might not be a foliation. Clearly
$F_0$ is a foliation. But $F_1$ is a foliation if and only 
if the projective connection on $M'$
satisfies a certain complicated condition on
its curvature, ensuring the existence of a suitably
large family of totally geodesic hypersurfaces, 
a local calculation which I leave to the reader.

On the other hand, if the projective connection on $M'$
has ``enough'' totally geodesic hypersurfaces, so that
$F_1$ is a foliation, then  each leaf  of $F_1$ 
projects to a immersed complex hypersurface in $M'$,
so that for every linear hyperplane $H \in \Proj{}T^*M'=M$,
there is a unique such complex hypersurface with tangent space $H$.
\end{example}

\begin{remark}
Suppose that $M'$ is a complex manifold with
holomorphic normal projective connection (see
Kobayashi and Nagano \cite{Kobayashi/Nagano:1964}
for the definition of normal).
We leave the reader to check that a
projective connection has ``enough'' totally geodesic hypersurfaces
(i.e. one through each point with each possible tangent
hyperplane, i.e. $F_1$ is a foliation) if and only
if the projective connection is flat.
\end{remark}

\begin{theorem}
Suppose that $M$ is a smooth complex projective
variety bearing a holomorphic double Legendre foliation.
Then $M=\Proj{}T^* M'$ for some smooth complex 
projective variety $M'$, and $M$ is the 
lift of a holomorphic projective connection on $M'$.
\end{theorem}
\begin{remark}
This theorem reduces the classification of double
Legendre foliations on smooth complex projective
varieties to that of holomorphic projective connections
satisfying the required curvature condition
to have ``enough'' totally geodesic hyperplanes.
\end{remark}
\begin{proof}
By theorem~\vref{theorem:Tabachnikov},
for any double Legendre foliation,
say with contact structure $\mathcal{I} \subset T^*M$,
there is a parabolic geometry $E \to M$ so that
$\mathcal{I} = E \times_P I$
for some $P$-module $I \subset \left(\LieG/\LieP\right)^*$. 
Since $P \subset G$ is a maximal parabolic
subgroup, by proposition~\vref{proposition:MaximalGivesSemicanonical},
$I$ is a semicanonical $P$-module.
Therefore the parabolic geometry on $M$
has a semicanonical module whose associated
Pfaffian system is not Frobenius. Apply
theorem~\vref{theorem:ContracanonicalDrop} to see that the
geometry must drop. Unless $M$ is isomorphic
to the model geometry, there is only
one space $A_n/P'$ that it can drop to, 
since there is only one parabolic subgroup
$P ' \subset A_n$ containing $P$, so $A_n/P'=\Proj{n}$
a projective connection on some $M'$. 
(To be precise, there are actually two such subgroups,
but there is only one up to outer automorphism.)
So $M$ is a lift of a projective connection on $M'$.
The lift of any projective connection to a parabolic
geometry modelled on the adjoint variety of $A_n$
is given in detail in example~\vref{example:AdjointAn},
and must be $M=\Proj{}T^* M'$. Since the
subbundle $F_1 \subset TM$ is a foliation,
the projective connection must have ``enough
hypersurfaces''.
\end{proof}

\section{Circles in parabolic geometries}

So far we have one method to force dropping
of Cartan geometries: semicanonical modules.
Next we will describe a different method to force dropping,
instead of using semicanonical modules. The method
of semicanonical modules works particularly well
on parabolic geometries modelled on $G/P$
with $P \subset G$ a maximal parabolic subgroup.
Our new method, the method of rational 
circles, will work only on the opposite
extreme: geometries modelled on $G/B$.
All parabolic geometries, with any model $G/P$, 
lift to geometries modelled on $G/B$, so it is natural
to focus on the $G/B$-geometries.

\begin{definition}
Suppose that $\alpha$ is a root of a complex
semisimple Lie group $G$. Let $\s{\alpha}$
be the Lie subalgebra of $\LieG$ generated
by the root spaces of $\alpha$ and $-\alpha$.
\end{definition}

\begin{remark}
Any root of any complex semisimple Lie algebra is reduced, so $\s{\alpha}$
is isomorphic to $\LieSL{2,\C{}}$ \cite{Serre:2001}.
\end{remark}

\begin{definition}
Suppose that $E \to M$ is a holomorphic
parabolic geometry, with Cartan connection
$\omega$, modelled on a generalized
flag variety $G/P$. Suppose that $\alpha$
is a positive simple root of $G/P$.
Define a Pfaffian system on $E$ by the equation
\[
\omega = 0 \pmod{\s{\alpha}}
\]
on tangent vectors. This Pfaffian system
has the fibers of $E \to E/B$ as Cauchy
characteristics, and is $B$-invariant, and therefore descends
to a unique Pfaffian system on $E/B$.
Call the maximal integral Riemann surfaces
of this system on $E/B$ \emph{$\alpha$-circles},
or just \emph{circles}.
\end{definition}

\begin{remark}
There is some danger of confusion here,
since $\omega$ is not actually defined on $E/B$,
and since the $\alpha$-circles of $M$ are
complex 1-dimensional submanifolds of $E/B$, not of $M$.
\end{remark}

\begin{remark}
Suppose that $G/P$ is a generalized flag variety and
that $G$ splits into a product
\[
G = G_1 \times G_2 \times \dots \times G_s
\]
of simple complex Lie groups. Then the Borel
subgroup $B \subset G$ has the form
\[
B = B_1 \times B_2 \times \dots \times B_s,
\]
where $B_i = B \cap G_i$. 
Let
\[
B'_i = P_1 \times P_2 \times P_{i-1} \times B_i \times P_{i+1} \times \dots \times P_s.
\]

If $\alpha$ is a positive
simple root of $\LieG_i$, then we can
define the circles by the equation
\[
\omega = 0 \pmod{\s{\alpha}}
\]
on $E$ as above, but we find that in fact
the fibers of $E \to E/B'_i$ are Cauchy 
characteristics for this linear Pfaffian
system. So in fact, we can define the
circles as Riemann surfaces on the various $E/B'_i$.
For our purposes in this paper, this observation
has no significance, but it should save computation
in examples.
\end{remark}

\begin{lemma}
In the model, $G/P$, all circles are rational.
\end{lemma}
\begin{proof}
In the model $G/P$, with the standard model
$G/P$-geometry, the $\alpha$-circles are precisely
the orbits in $G/B$ of the connected subgroup $\SL{2,\C{}}_{\alpha} \subset G$
whose Lie algebra is $\s{\alpha}$. 
Note that $\s{\alpha} \cap \LieB \subset \s{\alpha}$
is a Borel subalgebra. 
So the associated connected Lie subgroup $\SL{2,\C{}}_{\alpha}$
(which is either isomorphic to $\SL{2,\C{}}$ or to $\PSL{2,\C{}}$)
must act on the orbit as a complex semisimple Lie group
acting on a generalized flag variety. The orbit has one complex dimension.
Therefore the orbit is a rational curve. So 
in the model, all circles are rational.
\end{proof}

\begin{example}
Take $G/P=\Proj{2}$ and $G=\PSL{3,\C{}}$. 
There is one noncompact positive simple root
$\alpha$.
Note that $G/B=\Proj{}T\Proj{2}$.
In fact $\alpha$-circles are precisely the
lifts of complex projective lines in $\Proj{2}$
to the projectivized tangent bundle. We lift
each line $L=\Proj{1} \subset \Proj{2}$ 
by taking each point $p \in L$ to $T_p L \in \Proj{}T\Proj{2}=G/B$.
Clearly $\Proj{}T\Proj{2}$ is foliated by the lifts of lines.
\end{example}

\begin{remark}
For each positive simple root $\alpha$ of $G/P$,
the $\alpha$-circles foliate $E/B$. They are not
actually defined inside $M$, although each $\alpha$-circle
projects via a local biholomorphism to a 
Riemann surface in $M$. Therefore it is natural
to picture the $\alpha$-circles as (not necessarily compact)
curves in $M$.
\end{remark}

\begin{remark}
Another description of the $\alpha$-circles:
they are the leaves of the foliation 
$E \times_B \LieG_{-\alpha} \subset T(E/B) = E \times_B \left(\LieG/\LieB\right)$.
\end{remark}

\begin{theorem}[Brunella \cite{Brunella:2008}]\label{theorem:Brunella}
Suppose that $M$ is a compact K\"ahler manifold.
Suppose that $F \subset TM$ is a holomorphic foliation by 
(not necessarily compact) curves,
i.e. a rank 1 subbundle. Either (1) all of the leaves of $F$ are 
rational curves and $F^*$ is not pseudoeffective or (2)
none of the leaves of $F$ are rational and $F^*$ is pseudoeffective.
\end{theorem}

\begin{remark}
Brunella's theorem concerns holomorphic foliations with
singularities, but we will only consider nowhere singular
foliations, so case (2) above follows from Brunella's remarks
\cite{Brunella:2008} p. 55.
\end{remark}

\begin{proposition}\label{proposition:CircleDrops}
Suppose that $M$ is a compact K{\"a}hler manifold
bearing a holomorphic parabolic geometry modelled
on a generalized flag variety $G/P$.
Suppose that $B \subset P$ is a Borel subgroup.

Draw the Dynkin diagram of $P$, but
then change any cross (say corresponding to
some noncompact simple root $\alpha$)
to a dot if the $\alpha$-circles are rational.
In other words, change a cross to a dot 
just when the line 
bundle $E \times_B \LieG_{\alpha}$ on $E/B$
is not pseudoeffective.
Let $Q$ be the parabolic subgroup of $G$
whose Dynkin diagram we have just drawn.

Suppose $P' \subset G$ is a parabolic
subgroup containing $P$. Then $M$ drops
$M \to M'$ to a holomorphic parabolic geometry
modelled on $G/P'$ if and only if $Q \subset P'$.
\end{proposition}
\begin{proof}
By theorem~\vref{theorem:Brunella}, these line bundles
are pseudoeffective if and only if the $\alpha$-circles
are rational. 

If $M$ drops to $M \to M'$, then the
$\alpha$-circles lie inside the fibers of $M \to M'$,
i.e. inside generalized flag varieties $P'/P$. In these
generalized flag varieties, the induced $P'/P$-geometry
is the model geometry, and the $\alpha$-circles
are therefore rational curves. 

Conversely if
the $\alpha$-circles are rational curves, 
theorem~\vref{thm:CurvesForceDescent}
ensures that there is a drop $M \to M'$
so that all of the $\alpha$-circles
lie in the fibers of $M \to M'$. The fibers
are $P'/P$, some parabolic subgroup $P' \subset G$.
For $P'/P$ to contain all of the $\alpha$-circles, $P'$
must have $\alpha$ as a compact root.
\end{proof}

\begin{theorem}\label{theorem:BorelBracketClosed}
Suppose that 
\begin{enumerate}
\item $G$ is a complex semisimple Lie group with Borel subgroup $B \subset G$ and
\item $M$ is a compact K{\"a}hler manifold and
\item $E \to M$ is a holomorphic parabolic geometry modelled on $G/B$.
\end{enumerate}
Then  either (1) this parabolic geometry drops to some lower 
dimensional holomorphic parabolic geometry 
on a compact K{\"a}hler manifold or 
(2) for every $B$-submodule $I \subset \left(\LieG/\LieB\right)^*$,
the associated Pfaffian system $E \times_B I \subset T^*M$
is Frobenius.
\end{theorem}
\begin{proof}
By proposition~\vref{proposition:CircleDrops},
if the geometry does not drop, then for every
positive simple root $\alpha$, the line bundle $E \times_B \LieG_{\alpha}$
on $M$ is pseudoeffective. 
Every positive root $\alpha$ is a sum,
with nonnegative integer coefficients, of
positive simple roots. So for any positive root $\alpha$,
not necessarily simple, the line bundle $E \times_B \LieG_{\alpha}$
on $M$ is also pseudoeffective. 
Pick any $B$-submodule $I \subset \left(\LieG/\LieB\right)^*$.
Then
\[
\det I = \bigotimes_{\alpha} \LieG_{\alpha},
\]
where the tensor product is over positive roots $\alpha$ for
which $\LieG_{\alpha} \subset I$.
Therefore the line bundle
\(
E \times_B \det I
\)
on $M$ is pseudoeffective. By lemma~\vref{lemma:Demailly},
$E \times_B I \subset TM$ is Frobenius.
\end{proof}

\section{Example: second order scalar ordinary differential equations}

A path geometry is a geometric description of a system of 2nd order
ordinary differential equations.

\begin{example}
Take a 2nd order scalar order differential equation,
\[
\frac{d^2 y}{dx^2} = f\left(x,y,\frac{dy}{dx}\right).
\]
Pick a variable $p$, and consider the associated
foliation
\begin{align*}
dy &= p \, dx, \\
dp &= f\left(x,y,p\right) \, dx,
\end{align*}
whose leaves correspond to the solutions of the equation.
Also consider the foliation
\begin{align*}
dy &= 0, \\
dx &= 0,
\end{align*}
whose leaves correspond to the points $(x,y)$ of the
configuration space.
\end{example}

\begin{definition} A \emph{path geometry} on a complex manifold $M$
with $\dim_{\C{}} M = 3$ 
is a choice of 2 nowhere tangent holomorphic foliations on $M$ by
(not necessarily compact) curves, called \emph{integral curves}, and
\emph{stalks} respectively,
with both foliations being tangent to a (necessarily uniquely determined)
holomorphic contact plane field.
\end{definition}

\begin{remark}
In other words, a path geometry is a double Legendre foliation
of a 3-manifold.
\end{remark}

\begin{remark}
It turns out that near any point of $M$
there are local coordinates $x, \, y, \dot{y}$ on $M$ and
there is a holomorphic function $f\left(x,y,\dot{y}\right)$ for which the integral curves
are the solutions of
\[
dy = \dot{y} \, dx, \ d\dot{y} = f \, dx,
\]
while the stalks are the solutions of $dx=dy=0$.
Conversely, for any holomorphic function $f\left(x,y,\dot{y}\right)$,
these two holomorphic foliations are a path geometry.
\end{remark}

\begin{remark}
If we interchange the two foliations of a path geometry,
we obtain the \emph{dual path geometry}.
\end{remark}

\begin{theorem}[Cartan \cite{Cartan:1938}]
A holomorphic path geometry on a complex manifold $M$
determines and is determined by a holomorphic
parabolic geometry $E \to M$ modelled on the adjoint
variety $B_2/B$.  The holomorphic contact
structure is $E \times_P I$ for a semicanonical
$P$-submodule $I \subset \left(\LieG/\LieB\right)^*$.
\end{theorem}
\begin{remark}
See \cite{Bryant/Griffiths/Hsu:1995} for a detailed
exposition.
\end{remark}
\begin{remark}
We encountered this adjoint variety in example~\vref{example:adjoint}.
Suppose that $M$ is a complex manifold with path geometry,
and that $E \to M$ is the induced regular parabolic geometry
of Cartan's theorem. 
Labelling roots and $B$-modules as in example~\vref{example:An},
we can see that $E \times_B I_{12}$ is a holomorphic contact
structure. This is the contact structure of the path geometry.
The integral curves are circles of the root $\alpha_1$,
while the stalks are the circles of the root $\alpha_2$.
\end{remark}

Using our method of semicanonical modues, 
we can classify holomorphic path geometries on smooth complex 
projective varieties. By the method of circles, we find the complete classification on
compact K{\"a}hler manifolds.

\begin{theorem}\label{theorem:PathGeometriesDrop}
Suppose that $M^3$ is a connected compact K{\"a}hler 
manifold with a holomorphic path geometry. Then 
$M =\Proj{}TM'$, where $M'$ is a compact K{\"a}hler surface
with a holomorphic projective connection.
The stalks of $M \to M'$ are either the stalks or the integral curves of $M$.
The manifold $M'$ is 
\begin{enumerate}
\item
$\Proj{2}$ (and $M$ is the model $B_2/B$ with its
standard flat path geometry) or
\item
a complex surface with an unramified covering by the 
unit ball in $\C{2}$ 
(and $M$ is a quotient of an open set in the
model $B_2/B$, with its standard flat path geometry),
or
\item
a complex surface with an unramified holomorphic
covering by a 2-torus
(and the pullback projective connection on the torus is 
translation invariant).
\end{enumerate}
All of these possibilities for $M'$ occur. The
parabolic geometry on $M$ is the path
geometry associated to the geodesic 
equation of the projective connection on $M'$.
\end{theorem}
\begin{proof}
We have proven in a more general setting in section~\vref{section:DoubleLegendreFoliations}
that $M=\Proj{}TM'$, 
(keeping in mind that $\Proj{}TM' = \Proj{}T^*M'$).
The compact complex surfaces which bear projective
connections have been classified \cite{Kobayashi/Ochiai:1980}:
$M'=\Proj{2}$ or $M'$ is an unramified ball quotient or 
an unramified torus quotient.
The projective connections on these surfaces have also been classified
\cite{Klingler:1998,Klingler:2001} under the hypothesis
of local flatness. In case $M'=\Proj{2}$, the presence
of rational curves in $\Proj{2}$ ensures, by theorem~\vref{thm:CurvesForceDescent},
that the holomorphic projective connection on $M'$ is flat. 
In case $M'$ is covered by the
ball, Klinger's arguments in \cite{Klingler:2001} actually
go through without change, to prove local flatness.
Finally, if $M'$ is covered by a torus, then any projective connection
on $M'$ is translation invariant as shown in \cite{McKay:2008c}.
A local calculation (see Bryant, Griffiths and Hsu \cite{Bryant/Griffiths/Hsu:1995}) 
shows that the lift $M$ of a projective
connection on any complex surface $M'$ has
parabolic geometry given by the path geometry
of the geodesic equation of the projective connection.
\end{proof}

\begin{remark}
Another perspective: either the scalar second order 
ordinary differential equations that comprise the parabolic geometry
on $M$ are the geodesic equations of a holomorphic
projective connection on a complex surface, or else
they are the dual equations of such equations.
\end{remark}

\begin{remark}
The 2nd order ODE of the path geometry, in the 
case when $M'=\Proj{2}$ or a ball quotient, is locally equivalent to
\[
 \frac{d^2 y}{dx^2} = 0.
\]
If $M'$ is a surface covered by a torus, then the 
2nd order ODE is locally equivalent to
\[
 \frac{d^2 y}{dx^2} = p\left(\frac{dy}{dx}\right),
\]
for $p$ a polynomial of degree at most 3 with
constant coefficients (in linear holomorphic coordinates on 
the torus); see Cartan \cite{Cartan:1992} for proof. 
In either case, these equations are solvable by
quadratures.
\end{remark}

\section{Example: third order scalar ordinary differential equations}

Sato \& Yoshikawa \cite{Sato/Yoshikawa:1998} 
have results on third order ordinary 
differential equations similar to Cartan's above on
second order  ordinary 
differential equations. 
For any third order ordinary differential 
equation for one function of one
variable, say
\[
\frac{d^3y}{dx^3}=
f\left(x,y,\frac{dy}{dx},\frac{d^2y}{dx^2}\right),
\]
consider the manifold $M^4$ whose
coordinates are $x,y,p,q$, 
equipped with the exterior differential
system
\[
dy=p \, dx, \ dp = q \, dx, \ dq = f(x,y,p,q) \, dx.
\]
Sato and Yoshikawa put a parabolic geometry
on $M$. Their parabolic geometry is invariant
under ``contact transformations''. 

In this context, a contact transformation (in the sense
of Lie, not the sense of contact topology)
is any biholomorphism that preserves 
a certain complete flag of Pfaffian systems.
Let
\begin{align*}
\vartheta_1 &= dy-p \, dx, \\
\vartheta_2 &= dp-q \, dx, \\
\vartheta_3 &= dq-f(x,y,p,q) \, dx
\end{align*}
and let $\mathcal{I}_j$ ($j=1,2,3$) be the Pfaffian
system spanned locally by $\vartheta_1, \dots, \vartheta_j$.
Then a contact transformation in Lie's sense
is a local biholomorphism preserving all
of the Pfaffian systems $\mathcal{I}_j$.

\begin{figure}
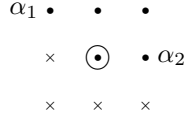

\centering{
\def\latticebody{%
\drop{} }
\[
\xy +(2,2)="o",0*\xybox{%
0;<1.5pc,0mm>:<0mm,1.5pc>::,{"o" \croplattice{-1}1{-1}1{-1}{1}{-1}1}
,"o"+(0,0)="a"*{\pin}*+!D{}
,"o"+(0,0)="b"*{\bigcirc}*+!L{}
,"o"+(1,0)="A"*{\pin}*+!D{} ,"o"+(0,1)="B"*{\pin}*+!L{}
,"o"+(1,1)="C"*{\pin}*+!L{} ,"o"+(-1,1)="D"*{\pin}*+!L{}
,"o"+(1,0)="A"*{\qquad \alpha_2}*+!D{} ,"o"+(-1,1)="B"*{\alpha_1 \qquad}*+!L{}
,"o"+(-1,0)="Ax"*{\pout}*+!D{} ,"o"+(0,-1)="Bx"*{\pout}*+!L{}
,"o"+(-1,-1)="Cx"*{\pout}*+!L{} ,"o"+(1,-1)="D"*{\pout}*+!L{} }
\endxy
\]
} \caption{The root spaces of the Borel subalgebra of
$B_2$.}\label{fig:latticeBTwo}
\end{figure}

The parabolic geometry of Sato and Yoshikawa is modelled on
$\Sp{4,\C{}}/B=C_2/B=B_2/B=\PO{5,\C{}}/B$, where $B$ is the
Borel subgroup. The root lattice of $B_2/B$ is drawn in figure~\ref{fig:latticeBTwo}. 
One can see the 4 positive roots, drawn as dots.
The 3 Pfaffian systems are associated to the sets of positive roots
\begin{align*}
& \alpha_1 \\
& \alpha_1, \alpha_1 + \alpha_2 \\
& \alpha_1, \alpha_1 + \alpha_2, \alpha_1 + 2 \, \alpha_2.
\end{align*}
The Dynkin diagram of the model is \(
 \begin{xy}
  0;/r.20pc/:
  (0,0)*+{\pout}="1";
  (5,0)*+{>}="2";
  (10,0)*+{\pout}="3";
  "1"; "3" **\dir{=};
  {\ar@{=}; "1";"3"};
  \end{xy}
\). 
Let's refer to a parabolic geometry with this model which arises
locally from a third order ordinary differential equation, following
the method of Sato and Yoshikawa, as a
\emph{third order ODE geometry}.
The fibers of the bundle map
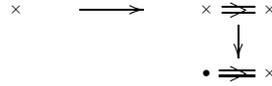
\[
\begin{xy}
  0;/r.20pc/: % Set the scale to be quite small
  (0,10)*+{\pout}="1";
  (10,10)*+{}="2";
  (20,10)*+{}="3";
  (30,10)*+{\pout}="4";
  (35,10)*+{>}="5";
  (40,10)*+{\pout}="6";
  (30,0)*+{\pin}="7";
  (35,0)*+{>}="8";
  (40,0)*+{\pout}="9";
  "2"; "3" **\dir{-};
  "4"; "6" **\dir{=};
  {\ar@{->}; "2";"3"};
  {\ar@{->}; "5"; "8"};
  "7"; "9" **\dir{=};
\end{xy}
\]
are the $e_3$-circles of the model (in the terminology of Sato \&
Yoshikawa) while those of
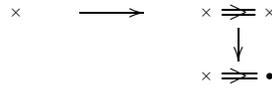
\[
\begin{xy}
  0;/r.20pc/: % Set the scale to be quite small
  (0,10)*+{\pout}="1";
  (10,10)*+{}="2";
  (20,10)*+{}="3";
  (30,10)*+{\pout}="4";
  (35,10)*+{>}="5";
  (40,10)*+{\pout}="6";
  (30,0)*+{\pout}="7";
  (35,0)*+{>}="8";
  (40,0)*+{\pin}="9";
  "2"; "3" **\dir{-};
  "4"; "6" **\dir{=};
  {\ar@{->}; "2";"3"};
  {\ar@{->}; "5"; "8"};
  "7"; "9" **\dir{=};
\end{xy}
\]
are the integral curves of the model (which Sato \& Yoshikawa call
$e_4$-circles).

\begin{example}
Suppose that $M'$ is a complex manifold,
of complex dimension 3,
with holomorphic parabolic geometry  $E \to M'$
modelled on the smooth quadric hypersurface $Q^3=\PO{5,\C{}}/P=B_2/P$.
For example, a holomorphic conformal structure
on $M'$ will impose such a holomorphic parabolic
geometry. Conversely, every parabolic geometry
modelled on $B_2/B$ 
imposes a holomorphic conformal structure, since
the group $P$ acts in the representation $\LieG/\LieP$
preserving a nondegenerate quadratic cone. 
Let $M$ be the set of all null
lines in the tangent spaces of $M'$. We can
easily see that $M=E/B$, $B \subset B_2$ the
Borel subgroup. Therefore $M$ is the lift of $M'$
to a $B_2/B$-geometry.
Moreover, if the parabolic geometry on $M'$ is
a holomorphic conformal structure, then the
parabolic geometry of the lift $M$ is precisely
the equation of circles in $M'$. We leave the 
reader to check these (purely local and elementary)
assertions of parabolic geometry.
\end{example}

\begin{remark}
The construction of a third order ordinary differential equation out
of a conformal structure has been well known since work of
W\"unschmann (see Chern \cite{Chern:1937,Chern:1940}, Dunajski \&
Tod \cite{Dunajski/Tod:2005}, Frittelli, Newman \& Nurowski
\cite{Frittelli/Newman/Nurowski:2003}, Sato \& Yoshikawa
\cite{Sato/Yoshikawa:1998}, Silva-Ortigoza \&
Garc{\'{\i}}a-God{\'{\i}}nez
\cite{SilvaOrtigoza/GarciaGodinez:2004}, W\"unschmann
\cite{Wunschmann:1905}). Identification of the local obstruction
to dropping with the Chern invariant is a long but straightforward
calculation (see Sato \& Yoshikawa \cite{Sato/Yoshikawa:1998}).
Hitchin \cite{Hitchin:1982} pointed out that a rational curve on a
surface with appropriate topological constraint on its normal bundle
must lie in a moduli space of rational curves constituting the
integral curves of a unique third order ordinary differential
equation with vanishing Chern invariant.
\end{remark}

\begin{theorem}\label{theorem:KahlerBTwoModB}
Suppose that $M$ is a compact K{\"a}hler 4-fold with
holomorphic parabolic geometry modelled on $B_2/B$.
Then the geometry on $M$ drops to 
a holomorphic parabolic geometry modelled
on
\begin{enumerate}
\item  the smooth quadric hypersurface $Q^3=B_2/P$ or
\item the projective space $\Proj{3}$ of null 2-planes in $\C{5}$,
\end{enumerate}
on a compact K{\"a}hler 3-fold $M'$.
In particular, if the parabolic geometry on $M$
is locally a third order ODE geometry, then 
this geometry is the equation of 
\begin{enumerate}
\item circles of a holomorphic conformal structure on $M'$ or
\item circles of a holomorphic Legendre connection on $M'$ 
\end{enumerate}
\end{theorem}
\begin{remark}
See Sato and Yoshikawa \cite{Sato/Yoshikawa:1998}
for the definition of Legendre connection.
\end{remark}
\begin{remark}
It is not known which compact K\"ahler 3-folds admit conformal
geometries, or admit Legendre connections.
\end{remark}
\begin{proof}
The model for third order ODE geometries is 
$B_2/B$, a quotient by a Borel subgroup,
so theorem~\vref{theorem:BorelBracketClosed} applies,
ensuring that the parabolic geometry drops
to one modelled on either
\(
 \begin{xy}
  0;/r.20pc/:
  (0,0)*+{\pout}="1";
  (5,0)*+{>}="2";
  (10,0)*+{\pin}="3";
  "1"; "3" **\dir{=};
  {\ar@{=}; "1";"3"};
  \end{xy}
\)
or
\(
 \begin{xy}
  0;/r.20pc/:
  (0,0)*+{\pin}="1";
  (5,0)*+{>}="2";
  (10,0)*+{\pout}="3";
  "1"; "3" **\dir{=};
  {\ar@{=}; "1";"3"};
  \end{xy}
\).
We call these $B_2/P_1$ and $B_2/P_2$ respectively.
The variety $B_2/P_1$ is the Dynkin diagram of the model
of a conformal geometry.
The variety $B_2/P_2$ is the Dynkin diagram of the model
of the parabolic geometry of a Legendre connection (again see Sato and Yoshikawa
\cite{Sato/Yoshikawa:1998}). We leave the
reader to check by examination of the  structure
equations of Sato and Yoshikawa \cite{Sato/Yoshikawa:1998} p. 1000 that 
if we take a parabolic geometry with either
of these two models on some complex manifold $M'$ 
and lift it, say to a complex manifold $M$, then
the lifted parabolic geometry on $M$ is the parabolic geometry associated
by Sato and Yoshikawa to the third order ODE of  the circles. 
\end{proof}

\begin{definition}
Recall that the \emph{Lie ball} is the noncompact
Hermitian symmetric space dual to the smooth quadric
hypersurface.
\end{definition}

\begin{theorem}
Suppose that $M$ is a smooth complex projective 4-fold with
a holomorphic third order ODE geometry.
Then $M$ is the set of null lines in the tangent
spaces of a 3-fold $M'$ with holomorphic conformal geometry.
The 3rd order ODE geometry on $M$ is the one
associated to the circles of $M'$.
The smooth complex projective 3-folds $M'$ which admit
holomorphic conformal structures are precisely
\begin{enumerate}
\item the quadric $Q^3$, with its standard flat conformal geometry,
\item 3-folds with unramified covering by an abelian 3-fold,
with any translation invariant conformal geometry,
\item 3-folds covered by the Lie ball with the standard
flat conformal geometry.
\end{enumerate}
\end{theorem}
\begin{proof}
Apply theorem~\vref{theorem:KahlerBTwoModB}
to ensure that the parabolic geometry on $M$ drops to a parabolic geometry
on some 3-fold $M'$. The parabolic geometry on $M'$ 
could be either a conformal structure or a Legendre connection.
Legendre connections admit a holomorphic
contact structure, of the form $E \times_{P_2} I$,
as we see from the
structure equations of Sato and Yoshikawa,
\cite{Sato/Yoshikawa:1998} p. 1000.
A contact structure is semicanonical and not Frobenius.  Therefore 
any Legendre connection on any smooth complex projective
variety must be isomorphic to
the model $\Proj{3}$ and this forces $M$ to be isomorphic
to its model, so drops to the model $Q^3$ of conformal
geometry.

Therefore we can assume that the parabolic geometry
on $M$ drops to a conformal geometry on a smooth
connected complex projective 3-fold $M'$.
The classification of smooth connected complex projective 3-folds
admitting conformal geometries is due to Jahnke and Radloff
\cite{Jahnke/Radloff:2004}. 
\end{proof}

\section{Conclusion}

We have demonstrated rigidity phenomena for a large class of
holomorphic geometric structures and holomorphic
exterior differential systems on smooth complex
projective varieties. Our motivation is the following conjecture.

\begin{conjecture}
Suppose that $G$ is a complex
simple Lie group and $P \subset G$ a maximal
parabolic subgroup. Suppose that $G/P$ is not
a compact Hermitian symmetric space (or, if 
$G/P$ is a compact Hermitian symmetric space, 
then suppose that $G$ is a proper subgroup of 
the identity component of the 
biholomorphism group of $G/P$).
Then up to isomorphism, the only holomorphic parabolic geometry 
modelled on $G/P$ on any compact K{\"a}hler manifold
is the standard flat parabolic geometry on $G/P$.
\end{conjecture}

More generally, one would like 
to construct explicitly all of the holomorphic Cartan
geometries on all compact K{\"a}hler manifolds.
The methods in this paper say nothing
about the parabolic geometries modelled on 
compact Hermitian symmetric spaces,
perhaps the most important type of
parabolic geometry \cite{Goncharov:1987,Klingler:2001}.

It might be possible to classify the semicanonical modules
of generalized flag varieties.  This is a complicated combinatorial
problem about root systems.

It is frustrating to have many results for smooth complex
projective varieties but so few for compact K{\"a}hler
manifolds. 
Any Cartan 2-plane field $\mathcal{V}$ on any 
5-dimensional complex manifold $M$ has a holomorphic 
quartic symmetric form on the
2-plane as an invariant; \cite{Cartan:30}. That quartic form 
has a discriminant, which is a holomorphic 
section of a positive power of the canonical bundle.
If the underlying 5-fold is compact K{\"a}hler, then
the canonical bundle is not pseudoeffective, as explained
above. Therefore no  positive power of the canonical bundle 
has any nonzero sections. So the discriminant vanishes, i.e.
the quartic has a multiple root at every point. 
The 2-plane field together with its brackets spans a 3-plane field. 
Similarly, Cartan defines a holomorphic quartic symmetric
form on the 3-plane field, which restricts to the
quartic on the two plane field. Again this quartic
can't have any nonvanishing invariants in classical invariant
theory, since these all occur in positive powers of 
the canonical bundle. By geometric invariant theory, the projectivized zero
locus of the quartic must therefore have a triple point or tacnode.
One might be able to find similar information about 
other invariants and thereby prove vanishing of curvature to prove that 
$G_2/P_1$ is the only compact K{\"a}hler 5-fold
bearing a holomorphic Cartan 2-plane field. It is already known
that the holomorphic Cartan 2-plane field on $G_2/P_1$ discovered by Cartan
is the only holomorphic Cartan 2-plane field on $G_2/P_1$ \cite{Biswas/McKay:2010}.

\bibliographystyle{amsplain}
\bibliography{rational-curves-parabolic}

\end{document}